\documentclass[preprint,12pt]{elsarticle}

\usepackage{amssymb}
\usepackage{amsmath}

\usepackage{graphicx} 
\usepackage{empheq}
\usepackage{algorithm}
\usepackage{algpseudocode}
\usepackage[table]{xcolor}
\usepackage{subcaption}
\usepackage{graphicx} 
\setcitestyle{authoryear, round}
\usepackage{enumitem}
\usepackage{hyperref}
\usepackage{array}
\usepackage{multirow} 
\usepackage{tikz}

\graphicspath{{figures/}}
\journal{}

\usetikzlibrary{shapes.geometric, arrows, fit}

\tikzstyle{startstop} = [rectangle, rounded corners, minimum width=3cm, minimum height=1cm, text centered, draw=black, fill=red!30]
\tikzstyle{process} = [rectangle, minimum width=4.5cm, minimum height=1cm, text centered, draw=black, fill=black!10]
\tikzstyle{arrow} = [thick,->,>=stealth]

\begin{document}

\begin{frontmatter}



\title{An enhanced heuristic framework for solving the Rank Pricing Problem}


\author[label1,label2]{Asunci\'on Jim\'enez-Cordero\corref{cor1}}
\ead{asuncionjc@uma.es}
\author[label1,label3]{Salvador Pineda}
\ead{spineda@uma.es}
\author[label1,label2]{Juan Miguel Morales}
\ead{juan.morales@uma.es}
\affiliation[label1]{organization={Group OASYS},
            addressline={Ada Byron Research Building, Arquitecto Francisco Peñalosa Street, 18}, 
            city={Málaga},
            postcode={29010}, 
            country={Spain}}
\affiliation[label2]{organization={Department of Mathematical Analysis, Statistics and Operations Research and Applied Mathematics (Area of Statistics and Operations Research), Universidad de Málaga},
            city={Málaga},
            country={Spain}}
\affiliation[label3]{organization={Department of Electric Engineering, Universidad de Málaga},
            city={Málaga},
            country={Spain}}
\cortext[cor1]{Corresponding author}
\begin{abstract}
The Rank Pricing Problem (RPP) is a challenging bilevel optimization problem with binary variables whose objective is to determine the optimal pricing strategy for a set of products to maximize the total benefit, given that customer preferences influence the price for each product. Traditional methods for solving RPP are based on exact approaches which may be computationally expensive. In contrast, this paper presents a novel heuristic approach that takes advantage of the structure of the problem to obtain good solutions. The proposed approach consists of two phases. Firstly, a standard heuristic is applied to get a pricing strategy. In our case, we choose to use the Variable Neighborhood Search (VNS), and the genetic algorithm. Both methodologies are very popular for their effectiveness in solving combinatorial optimization problems. The solution obtained after running these algorithms is improved in a second phase, where four different local searches are applied. Such local searches use the information of the RPP to get better solutions, that is, there is no need to solve new optimization problems. Even though our methodology does not have optimality guarantees, our computational experiments show that it outperforms Mixed Integer Program solvers regarding solution quality and computational burden.
\end{abstract}




\end{frontmatter}



\section{Introduction} \label{sec: introduction}

In recent times, pricing optimization problems have attracted attention within the Operations Research community, \citep{calvete2019rank, labbe2013bilevel}. Essentially, these models have to determine the optimal prices of a series of available products to maximize the benefits of a company at the same time that the customer preferences and their budget values are taken into account. In essence, understanding what customers prefer and can afford, or in other words, maximizing customer satisfaction within their budget limits, is crucial when setting the optimal prices for products.

This problem, which \emph{a priori} may seem easy, is very difficult to solve. Setting high prices may increase revenue, but if customers cannot afford the products or their preferences are not met, sales may drop, leading to low or zero profit. Conversely, setting low prices ensures purchases but limits revenue. Therefore, balancing these opposing factors to maximize profit is a difficult challenge.

To solve such an issue, in this paper, we focus on the so-called Rank Pricing Problem (RPP). This model has been proven to be NP-complete in \citep{rusmevichientong2006nonparametric}. The Rank Pricing Problem, as formulated here, was first introduced by \citep{rusmevichientong2006nonparametric}. RPP is a combinatorial optimization model that can be formulated as a bilevel optimization problem, \citep{dempe2020bilevel}. Given the fact that each customer has a fixed budget and a preference value for each product, the main objective of the RPP is to decide the price for the products in such a way that the total benefit of the company is maximized. The upper-level problem aims to maximize the revenue. The lower-level models (one per customer) have to maximize the satisfaction of the incumbent customer.

While the upper-level problem solves a linear optimization model (once the decision variables of the lower level are fixed), the lower-level formulation has to solve an integer linear optimization problem with binary variables. Such integer nature, together with the intrinsic nonlinearity of bilevel problems, increases the difficulty of finding the optimal solution. 

\section{Literature review}

Several exact techniques have been proposed in the literature to solve the Rank Pricing Problem. We can highlight, for instance, the works in \citep{calvete2019rank} and \citep{dominguez2021rank}. In these papers, a single-level formulation and different linearizations are proposed. Some valid inequalities are introduced to tighten the linear relaxation of the model and consequently improve the resolution of the problem. In addition, the works in \citep{calvete2024novel, dominguez2021rankpricing} deal with a variant of the RPP that considers ties in the preferences. Here some valid inequalities are also introduced. Similar optimization problems have been handled with exact algorithms. This is the case of the work in \citep{bucarey2021models} where a Benders decomposition with some valid inequalities is proposed in a related pricing model. A Benders decomposition algorithm was also applied in \citep{bertsimas2019exact}, where an optimization problem with ranking-based customer preferences is considered. Even though previous works based on exact approaches have presented outstanding results, it may be difficult and very time-consuming to find the right valid inequalities necessary to find a (possible) optimal solution.


Regarding the heuristic approaches, the work in \citep{aggarwal2004algorithms} proposes an algorithm for solving a modification of the RPP. Recently, a heuristic approach was proposed in \citep{calvete2024evolutionary} to solve the RPP. A genetic algorithm combined with a local search is applied to refine incumbent pricing strategies. Two key aspects of this approach are highlighted below. First, the initialization process is not purely random. Instead, the first price vector is constructed using a greedy strategy. Second, they implement a local search that runs immediately after the genetic algorithm. More precisely, for a given price vector, the price of each product is iteratively adjusted to match the price of any other available product. The updated price vector is then used to solve the lower-level problems and evaluate the upper-level objective function. The vector yielding the highest revenue is retained. While this method produces good solutions, it requires solving a large number of optimization problems.

In this paper, we improve upon the work by \citep{calvete2024evolutionary} and introduce a novel heuristic approach that leverages the inherent structure of the Rank Pricing Problem (RPP) to obtain high-quality solutions efficiently. Our strategy consists of two main steps. First, a standard heuristic algorithm from the literature is executed. Second, the resulting solution is refined using four new local search techniques. Unlike the approach in \citep{calvete2024evolutionary}, our improvements do not require solving additional optimization problems, thereby reducing computational overhead. As our numerical results demonstrate, the proposed heuristic method frequently produces solutions close to optimal, even though global optimality cannot be guaranteed. Consequently, as with most heuristic methods, we can only ensure local optimality.

For the base heuristic algorithm, we consider two well-established techniques: Variable Neighborhood Search (VNS) \citep{mladenovic1997variable, hansen2006variable, lazic2010variable} and genetic algorithms \citep{holland1975adaptation, kramer2017genetic, lambora2019genetic}. Since their introduction, both methods have been widely applied to generic mixed-integer problems \citep{rothberg2007evolutionary, deep2009real} and combinatorial optimization problems, including routing \citep{park2021waiting, xu2018variable}, healthcare management \citep{lan2021survey, tahir2022novel}, scheduling \citep{silva2023solving, squires2022novel}, and location problems \citep{cazzaro2022variable, jaramillo2002use}. While the genetic algorithm has already been used to solve the Rank Pricing Problem in \citep{calvete2024evolutionary}, to the best of our knowledge, VNS has not yet been applied to this problem.

The remainder of the paper is structured as follows: Section \ref{sec: The Rank Pricing Problem} describes the main aspects of the Rank Pricing Problem. Section
\ref{sec: methodology} details the different strategies used to solve this problem, including our proposal. Section \ref{sec: numerical experience} is devoted to computational experiments. Finally, we finish with some conclusions and further research in Section \ref{sec: conclusions}.

\section{The Rank Pricing Problem}\label{sec: The Rank Pricing Problem}

In this paper, we focus on a pricing model based on ranked preferences, designed for unit-demand customers with a positive budget and who have access to an unlimited supply of products. This version of the RPP can be described as follows: let $\mathcal{I}= \{1, \ldots, I\}$ be a set of products, and let $\mathcal{K}= \{1, \ldots, K\}$ denotes the set of customers. In addition, each customer $k\in \mathcal{K}$ has an available budget, $b^k$, and a preference value for each product $i\in \mathcal{I}$, denoted by $s_i^k$, \textcolor{black}{both of which are assumed to be known}. In this way, for a fixed customer $k$, $s_{i_1}^k>s_{i_2}^k$ implies that customer $k$ prefers to buy product $i_1$ rather than product $i_2$. Given this information, one can formulate RPP as the following bilevel optimization problem:

  \begin{subequations} \label{eq: upper level}
  \begin{empheq}[left=\empheqlbrace]{align}
    \max\limits_{\boldsymbol{p}}&\sum\limits_{k\in \mathcal{K}}\sum\limits_{i\in \mathcal{I}}p_ix_i^k  \label{eq: upper level objective function}\\
    \text{s.t. }  & p_i\geq 0,\, i \in \mathcal{I}
    \end{empheq}
  \end{subequations}

  where for all $k\in \mathcal{K}$, $x_i^k$ is the optimal solution of the problem:

     \begin{subequations} \label{eq: lower level}
  \begin{empheq}[left=\empheqlbrace]{align}
    \max\limits_{\boldsymbol{x^k}} & \sum\limits_{i\in \mathcal{I}}s_i^kx_i^k \label{eq: lower level objective function}  \\
    \text{s.t. }  & \sum\limits_{i\in \mathcal{I}} x_i^k\leq 1 \label{eq: lower level first constraint} \\
    &\sum\limits_{i\in \mathcal{I}} p_ix_i^k\leq b^k \label{eq: lower level second constraint}\\
    &x_i^k \in\{0, 1\},\,  i \in \mathcal{I}.
    \end{empheq}
  \end{subequations}

The vector of product prices, $\boldsymbol{p}= (p_1, \ldots, p_I)^\intercal$, forms the upper-level decision variables where $p_i$ takes non-negative values for all elements $i$ in the set $\mathcal{I}$. For each customer, $k\in\mathcal{K}$, the lower-level decision variables, $\boldsymbol{x^k} = (x_1^k, \ldots, x_I^k)^\intercal$ are binary variables that indicate whether customer $k$ chooses to buy product $i$ or not. To simplify notation, we denote by $\boldsymbol{x} = (\boldsymbol{x^1}|\ldots|\boldsymbol{x^K})$ the $I\times K$ matrix that includes the decision variables for the $K$ lower-level problems. We also denote by $\phi(\boldsymbol{p}, \boldsymbol{x})$ the objective function \eqref{eq: upper level objective function} evaluated at vector $\boldsymbol{p}$ and matrix $\boldsymbol{x}$. In addition, the lower-level problem \eqref{eq: lower level} for customer $k$ and a fixed vector of prices $\boldsymbol{p}$ is denoted as $LL^k(\boldsymbol{p})$. In the cases where we want to highlight that the optimal solution of $LL^k(\boldsymbol{p})$ is obtained for a given vector $\boldsymbol{p}$ we use $\boldsymbol{x^k}(\boldsymbol{p})$. 

The objective of the upper-level problem \eqref{eq: upper level} is to maximize the total benefit of the company for given choices of the customers, that is, for given values of the lower-level decision variables $\boldsymbol{x}$. On the other hand, the lower-level problem has to decide which product each customer purchases so that global satisfaction is maximized, as can be seen in \eqref{eq: lower level objective function}. This decision should be taken following some constraints. Firstly, constraint \eqref{eq: lower level first constraint} indicates that each customer can buy at most one product. Then, for a given vector of prices $\boldsymbol{p}$, constraint \eqref{eq: lower level second constraint} indicates that customer $k$ can buy product $i$, i.e., $x_i^k=1$, if and only if the price chosen for product $i$, $p_i$ is smaller or equal than the budget of customer $k$, $b^k$.

It has been shown (see \citep{dominguez2021rank} for more details) that for a given vector of upper-level decision variables, $\boldsymbol{p}$, there exists a unique optimal solution of the lower-level problem, $\boldsymbol{x^k},\, k\in \mathcal{K}$. On the other hand, with the lower-level decision variables $x_i^k$, $i\in \mathcal{I}, k\in \mathcal{K}$, fixed, we can observe that the optimal prices $p_i$ belong to the set of possible budgets $b^k,\, k\in \mathcal{K}$. This is due to the fact that if a customer $k$ chooses to buy product $i$, i.e., $x_i^k=1$, then the price of the product is upper bounded by the budget $b^k$ because of constraint \eqref{eq: lower level second constraint}. In addition, the objective function of the upper-level problem \eqref{eq: upper level objective function} aims to get the maximum benefit. Hence, the upper-level decision maker wants to set the maximum price so that the customers can afford this product. Such maximum price value, then, has to coincide with one of the budgets. Otherwise, the objective function does not reach its maximum. This statement, which has already been noted in \citep{dominguez2021rank} and \citep{rusmevichientong2006nonparametric}, is key to our approach because even though \emph{a priori}, the feasible upper-level decision variables take continuous values, the optimal solution is attained in the discrete grid provided by the available budgets.

\section{Methodology}\label{sec: methodology}

In the following sections, we analyze the different approaches we use throughout the paper. Section \ref{subs: exact approach} outlines how the RPP problem is solved in an exact manner, whereas Section \ref{subs: naive} details a naive approach for this purpose. Then, Section \ref{subs: base algorithms} enumerates the different base algorithms used for the heuristic proposal and Section \ref{subs: initialization procedure} describes the initialization procedures. Finally, Section \ref{subs: local improvements} details the local searches used to improve the incumbent solutions.

\subsection{Exact approach} \label{subs: exact approach}

The first methodology we present here is the so-called \emph{exact} approach. As its name indicates, using this strategy, we search for the optimal solution of the bilevel optimization problem \eqref{eq: upper level} - \eqref{eq: lower level} using an off-the-shelf solver. To this aim, we utilize the single-level reformulation extracted from \citep{dominguez2021rank} and that can be seen below:

\begin{subequations} \label{eq: single-level reformulation}
  \begin{empheq}[left=\empheqlbrace]{align}
    \max\limits_{\boldsymbol{v},\, \boldsymbol{x}} & \sum\limits_{k\in \mathcal{K}}\sum\limits_{i\in \mathcal{I}} \left(\sum\limits_{m\in \mathcal{M}}b^m v_i^m\right)x_i^k \label{eq: single level objective function} \\
    \text{s.t. }  & \sum\limits_{m\in M} v_i^m\leq 1, \forall i\in I \label{eq: single level first constraint}\\
    & \sum\limits_{i\in I} x_i^k\leq 1, \forall k \in K  \label{eq: single level second constraint}\\
    &x_i^k\leq \sum\limits_{m\in M} v_i^m, \forall k\in K, i\in I. \label{eq: single level third constraint}\\
    &\sum\limits_{j\in I} s_j^kx_j^k\geq s_i^k\sum\limits_{m\in M} v_i ^m, \forall k\in K, i\in I \label{eq: single level fourth constraint}\\
    &v_i^m, x_i^k \in\{0, 1\}, \forall k\in K, i\in I, m\in M. \label{eq: single level fifth constraint}
    \end{empheq}
  \end{subequations}

Formulation \eqref{eq: single-level reformulation} is based on the fact that the optimal prices are attained in the budget values. Indeed, since two different customers can have the same budget value, we need to define a new set, $\mathcal{M} = \{1, \ldots, M\}$, that contains the indexes of the unique values of the sorted budgets, where by \emph{unique} we mean that values that are repeated are included only once. Therefore, we have that $b^{m_1}< b^{m_2}$ if $m_1<m_2$ for $m_1, m_2 \in \mathcal{M}$, and the set $\mathcal{B} = \{b^m, m\in \mathcal{M}\}$. Using this notation, one can define new binary variables $v_i^m$ in such a way that $v_i^m = 1$ if and only if the price $p_i$ takes the value of the budget $b^m$. This way, we can say that $p_i = \sum\limits_{m\in \mathcal{M}} b^m v_i^m$. For notation purposes, we denote by $\boldsymbol{v}$ as the $I\times M$ matrix whose entry $(i, m)$ is the decision variable $v_i^m$. 

This model is an integer nonlinear optimization problem with binary variables $\boldsymbol{v}$ and $\boldsymbol{x}$. The nonlinearity of the problem is due to the quadratic term in the objective function coming for the product of $v_i^m$ and $x_i^k$. The objective function \eqref{eq: single level objective function} is equivalent to \eqref{eq: upper level objective function}. In other words, the aim is to maximize the total revenue of the company. Constraint \eqref{eq: single level first constraint} indicates that each product $i\in \mathcal{I}$ must have only one price, then at most one binary variable $v_i^m$ can take the value 1. Then, constraint \eqref{eq: single level second constraint} is the same as \eqref{eq: lower level first constraint}. It indicates that one customer can choose at most one product. In addition, thanks to constraint \eqref{eq: single level third constraint}, we can see that if customer $k$ buys product $i$, i.e., $x_i^k = 1$, then $\sum\limits_{m\in \mathcal{M}}v_i^m= 1$, and consequently, customer $k$ can afford product $i$, and its value is set to one of the available budgets in $\mathcal{B}$. Equivalently,  if there exists a product $i$ whose price $p_i$ cannot be set to any of the values in $\mathcal{B}$, then $v_i^m = 0, \forall m \in \mathcal{M}$, and none of the customers $k$ can buy this product, that is to say, $x_i^k = 0, \forall k \in \mathcal{K}$.  Next, constraint \eqref{eq: single level fourth constraint} states that if customer $k$ can afford product $i$ (whose price is set to any of the budgets $b^m$), then customer $k$ can buy a product $j$ that she prefers the same or more than $i$. In other words, this constraint assures that customer $k$ buys the product he likes the most whenever he can afford it. Finally, constraint \eqref{eq: single level fifth constraint} indicates the binary character of the decision variables. Problem \eqref{eq: single-level reformulation} will be solved using an off-the-shelf solver. 

\subsection{Naive approach} \label{subs: naive}

Since the optimal solution of the upper-level decision variables, $p_i,\, i \in \mathcal{I}$ is attained at the set of available budgets, $\mathcal{B}$, then one possible strategy to obtain such an optimal solution is to carry out a simple exploration of the feasible solutions. In particular, this strategy works as follows: while a stopping criterion is not reached, randomly sample different vectors of product prices $\boldsymbol{p}^{\ell}$. Each price value $p_i^{\ell}, i\in \mathcal{I}$ is independently taken from the available budget set $\mathcal{B}$. Then, for each fixed vector, $\boldsymbol{p}^{\ell}$, we solve the lower-level problems $LL^k(\boldsymbol{p}^{\ell})$ for all the customers $k\in \mathcal{K}$. The optimal decision variables of these problems, $\boldsymbol{x^k} (\boldsymbol{p}^{\ell}), \forall k \in \mathcal{K}$, as well as the fixed vector of prices $\boldsymbol{p}^{\ell}$ will be used to evaluate the upper-level objective function $\phi(\boldsymbol{p}^{\ell},\, \boldsymbol{x^k} (\boldsymbol{p}^{\ell}))$. Assume that when the stopping criterion is met, we have a set $\mathcal{L}$ of vectors of the form $\boldsymbol{p}^{\ell}$. Finally, we consider as the optimal objective value, $\phi^{max}$, the highest value $\phi(\boldsymbol{p}^{\ell},\, \boldsymbol{x^k} (\boldsymbol{p}^{\ell}))$ where $\boldsymbol{p}^{\ell}\in\mathcal{L}$. In other words, $\phi^{max}:= \max\limits_{\boldsymbol{p}^{\ell}\in\mathcal{L}}\phi(\boldsymbol{p}^{\ell},\, \boldsymbol{x^k} (\boldsymbol{p}^{\ell}))$ and $\boldsymbol{p}^{max}$ is the argument that maximizes these values. Algorithm \ref{alg: naive} shows a pseudocode of this approach.

\begin{algorithm}
\caption{\emph{Naive}}\label{alg: naive}
\begin{algorithmic}
\State \textbf{Input:} stopping criterion, $\mathcal{B}$, $\mathcal{I}$, and $\mathcal{K}$.

\While{stopping criterion is not reached}
\begin{enumerate}[label={\arabic*)}]
    \item Randomly generate vectors $\boldsymbol{p}^{\ell}$ from the set of budgets $\mathcal{B}$.
    \item For each customer $k$, and for the fixed prices $\boldsymbol{p}^{\ell}$, solve $LL^k(\boldsymbol{p}^{\ell})$ and obtain the optimal solutions $\boldsymbol{x}^k(\boldsymbol{p}^{\ell})$.
    \item Compute $\phi\left(\boldsymbol{p}^{\ell},\, \boldsymbol{x^k} (\boldsymbol{p}^{\ell})\right)$.
\end{enumerate}
\EndWhile
\vspace*{0.1cm}
 \State 4) The set of all the vectors of prices generated is denoted by $\mathcal{L}$.\\
    5) Set $\phi^{max} :=  \max\limits_{\boldsymbol{p}^{\ell}\in\mathcal{L}}\phi(\boldsymbol{p}^{\ell},\, \boldsymbol{x^k} (\boldsymbol{p}^{\ell}))$.\\
    6) Denote by $\boldsymbol{p}^{max}$ the argument that maximizes the previous expression.
\vspace*{0.1cm}
\State \textbf{Output:} $(\boldsymbol{p}^{max}, \phi^{max})$.
\end{algorithmic}
\end{algorithm}

This naive strategy is very simple. It essentially reduces to explore different upper-level feasible solutions and solve one lower-level problem per customer. Nevertheless, finding a solution close to the optimal may be intractable even for small-size instances. 

{
\renewcommand\thesection{\color{black}\arabic{section}}

\subsection{Base algorithms} \label{subs: base algorithms}

In this paper, we opt for two well-known heuristics, namely the VNS and the genetic algorithms (this one already applied to the RPP in \citep{calvete2024evolutionary}). These strategies are appropriate for solving the RPP given its combinatorial nature and the fact that the optimal price decisions are attained in the budgets. The two base algorithms we choose in this paper are based on the concept of \emph{neighborhood}. However, each of them considers this definition in a different manner. Sections \ref{subsubs: VNS approach} and \ref{subsubs: genetic approach} analyze both strategies in the context of the RPP and study their differences and similarities in terms of the neighborhood definition.

Let us first introduce an illustrative example that will be used in the next sections. Assume given an instance of the problem \eqref{eq: upper level}-\eqref{eq: lower level} with $K = 8$ customers and $I = 2$ products. Table \ref{tab: data illustrative example} shows the budgets for each customer and the preference matrix. We observe, for instance, that if we focus on the first customer, $k = 1$, he prefers product 1 to product 2 because $s_1^1 = 2$ is greater than $s_2^1= 1$. Furthermore, he can afford a product whose price is less or equal to $b^1 = 18$. Therefore, in this particular case, the set of available budgets is $\mathcal{B} = \{18, 27, 34, 42, 50, 66\}$. 

\begin{table}[htb!]
    \centering
    \begin{tabular}{c|cccccccc}
    \hline
        $k$ & 1 & 2 & 3 & 4 & 5 & 6 & 7 & 8 \\
        \hline
        $s_1^k$ & 2 & 2 & 1 & 1 & 1 & 2 & 2 & 1 \\
        
        $s_2^k$ & 1 & 1 & 2 & 2 & 2 & 1 & 1 & 2 \\
        \hline
        $b^k$ & 18 & 66 & 27 & 34 & 66 & 50 & 42 &42 \\
        \hline
    \end{tabular}
    \caption{Data for the illustrative example.}
    \label{tab: data illustrative example}
\end{table}

Figure \ref{fig: feasible solutions illustrative example} shows the set of all the feasible solutions for this instance. The points marked with stars indicate the optimal solution. In this case, we have multiple solutions that correspond to the set $\{(34, 66), (50, 34), (66, 34)\}$. The color of the different points indicates the objective value. It has been grouped in different colors. For instance, the objective values between 130 and 152 are colored in black, whereas the values between 218 and 240 are highlighted in red. In addition, within each color range, the darker the color, the higher the objective values. For example, in the black range, grey colors are close to a value of 130, whereas dark black is associated with a value of 152.

\begin{figure}[htb!]
    \centering
    \includegraphics[scale = 0.6]{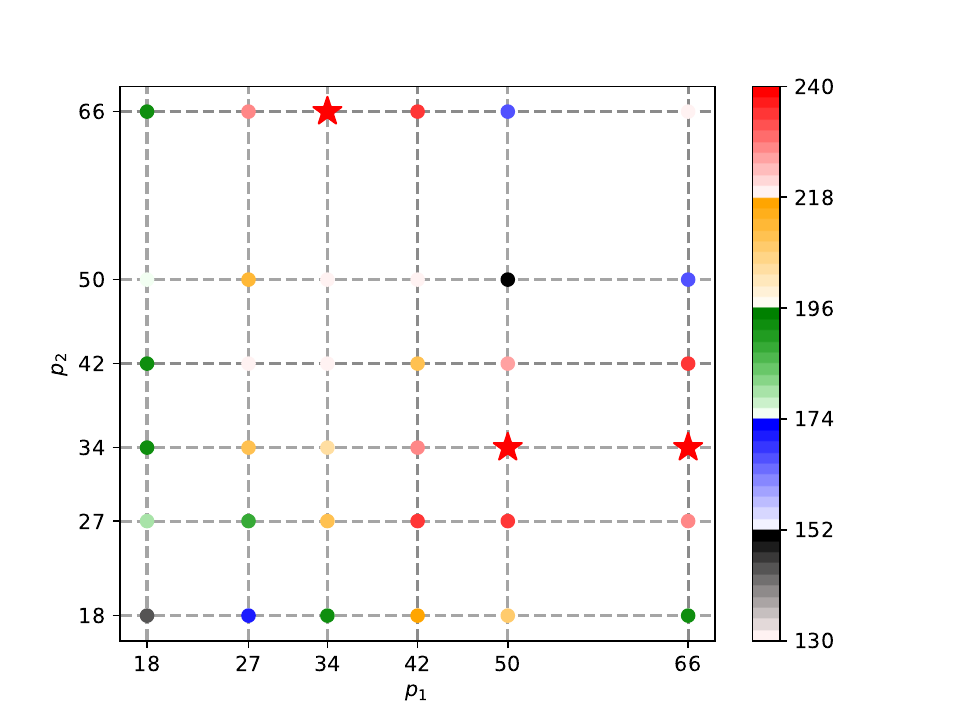}
    \caption{Feasible solutions for the illustrative example.}
    \label{fig: feasible solutions illustrative example}
\end{figure}


\subsubsection{The VNS-based approach} \label{subsubs: VNS approach}
In this section, we  explain how the VNS methodology is used for solving the RPP. Recall that the set $\mathcal{B}$ is formed by the set of unique sorted budgets. Indeed, one may think that changing the price of a product from $b^m$ to $b^{m+1}$ may have similar objective values. However, this is not always the case because a small change in the price of one product may cause big changes in the objective function.

Let us analyze this fact with the illustrative example of Figure \ref{fig: feasible solutions illustrative example}. We focus, for instance, on the vector of prices located at (50, 50). This point has an objective value between 130 and 152 (due to its black color). We observe that if we move one position in the upward direction, going to point (50, 66), this new point has a similar objective value, even though it is in a different color range. However, this is not the case if we move downwards to the point (50, 42) where we obtain a much larger objective value (red color). Indeed, with this movement, we move to a much better feasible solution in terms of the objective value. This effect is one of the consequences of the discrete nature of the RPP.

Now, we introduce how the basic version of the VNS approach works. Given an initial point, the VNS has two phases, the local-improvement phase and the perturbation phase. In the first one, a local optimum is found starting from the initial point. Then, we perturb the local optimum in a given neighborhood to have a new point and repeat the process. This process is repeated until a stopping criterion is met. It is, therefore, essential to define the neighborhoods considered in our proposal. Assume given a vector of prices $\boldsymbol{p} = (b^{m_1}, \ldots, b^{m_I})$. Note that $b^{m_i}\in \mathcal{B}, i \in \mathcal{I}$ denotes the budget value fixed for the price of product $i$. In addition, we want to highlight that two different products may have the same budget value; that is to say, it may happen that $m_i = m_j,$ with $i, j \in \mathcal{I}$ and $i\neq j$. As a consequence, $b^{m_i} = b^{m_j}$. We define the neighborhood around $\boldsymbol{p}$ of radius $r$, $r = 1, \ldots, R$, $N_r(\boldsymbol{p})$, as the following set:

\begin{footnotesize}
\begin{equation}\label{eq: neighborhood definition}
    N_r(\boldsymbol{p})= \left\{ (b^{n_1}, \ldots, b^{n_I})\in\mathcal{B}^I: n_i = m_i - r, \ldots, m_i -1,\, m_i,\, m_i+1, \ldots, m_i +r, \forall i \in \mathcal{I}\right\}.
\end{equation}
\end{footnotesize}

Returning to the illustrative example, we show in Figure \ref{fig: neighborhoods illustrative example} the neighborhoods around point $(50, 50)$ of radius $r= 1$ and $r=2$. In particular, the points that belong to $N_1((50, 50))$ are located in the black solid-line rectangle, whereas the points of the neighborhood $N_2((50, 50))$ are represented with a black dashed-line rectangle. 

\begin{figure}[htb!]
    \centering
    \includegraphics[scale = 0.7]{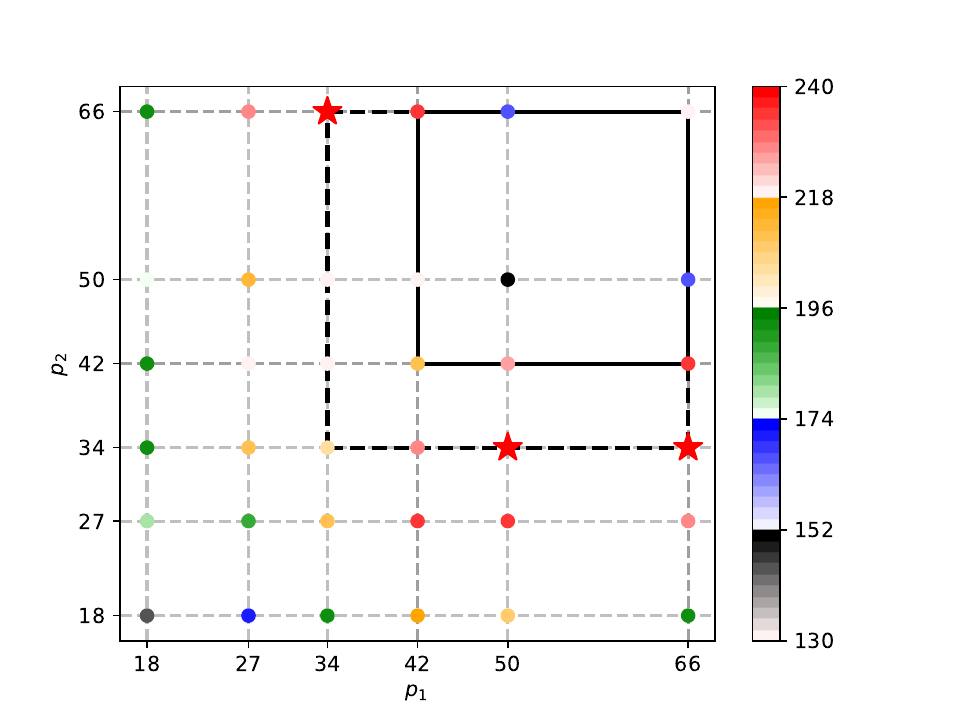}
    \caption{$N_1((50, 50)$ and $N_2((50, 50)$ for the illustrative example.}
    \label{fig: neighborhoods illustrative example}
\end{figure}

Once the definition of the neighborhoods is given, we will describe how our VNS-based proposal works: we first assume given an initial set $\mathcal{L}$ of $L_0$ vectors of prices, $\mathcal{L} = \{\boldsymbol{p}^\ell,\, \ell = 1, \ldots, L_0\}$. The way in which this initial population is generated is discussed in Section \ref{subs: initialization procedure}. Then, in the local-improvement phase, for each fixed vector $\boldsymbol{p}^\ell \in \mathcal{L}$, we solve the lower-level problems $LL^k(\boldsymbol{p}^\ell),\, k\in\mathcal{K}$, get the optimal solutions $\boldsymbol{x}^k(\boldsymbol{p}^\ell),\, k\in\mathcal{K}$, and evaluate the upper-level objective function $\phi(\boldsymbol{p}^\ell, \boldsymbol{x}^k(\boldsymbol{p}^\ell))$ that is denoted by $\phi^\ell, \ell = 1, \ldots, L_0$. Let us denote by $\phi^{max}$ the maximum objective value so far found. That is to say, $\phi^{max} = \max_{\ell = 1, \ldots, L_0} \phi^\ell$. In the following step, the $Q$ vectors of prices with the highest objective values in $\mathcal{L}$ are selected. We denote by $\mathcal{Q}$ the set that contains such vectors of prices, i.e., $\mathcal{Q} = \left\{ \boldsymbol{p}^{\ell_q}, q = 1, \ldots, Q: \phi^{\ell_1} \geq \ldots \geq \phi^{\ell_Q} \right\}$. Note that $\mathcal{Q}$ is a subset of $\mathcal{L}$. Next, in the perturbation phase, we randomly select $T$ points in the neighborhoods of radius $r = 1$ around the prices in $\mathcal{Q}$. In other words, we select $T$ points belonging to $N_1(\boldsymbol{p})$ for $\boldsymbol{p}\in \mathcal{Q}$. Let us denote by $\mathcal{T}$ the new set of $T$ price vectors. We update the set $\mathcal{L}$ by including the new points, i.e., $\mathcal{L}: = \mathcal{L} \cup \mathcal{T}$.  Later, we solve the lower-level problems $LL^k(\boldsymbol{p^t})$ for $\boldsymbol{p^t}\in \mathcal{T}$, and with the optimal decision variables so-obtained, $\boldsymbol{x}^k(\boldsymbol{p^t})$, we evaluate $\phi^t:= \phi(\boldsymbol{p^t}, \boldsymbol{x}^k(\boldsymbol{p^t}))$. If there is an objective value among these new points in $\mathcal{T}$ with a higher value than $\phi^{max}$,  or equivalently, if there exists an index $t$ in such a way that $\phi^t>\phi^{max}$, then we replace $\phi^{max}$ by $\phi^t$ and update the set $\mathcal{Q}$ by searching in the new set $\mathcal{L}$ the $Q$ vector of prices with the highest objective values. Otherwise, if $\phi^{max}\geq \phi^t, \forall t = 1, \ldots, T$, then we also update the set $\mathcal{Q}$ and increase in one unit the radius of the neighborhood to explore new points. Equivalently, we set $r:=r+1$ and generate new $T$ vectors of prices around the points belonging to the updated set $\mathcal{Q}$. This process is repeated until a stopping criterion is reached. At the end of the algorithm, we consider as the optimal objective value the value $\phi^{max}$ found so far. In addition, the optimal solution associated with such objective value is denoted by $\boldsymbol{p}^{max}$. For a detailed explanation of our approach, see Algorithm \ref{alg: VNS}.

\begin{algorithm}
\caption{\emph{VNS}}\label{alg: VNS}
\begin{algorithmic}
\State \textbf{Input:} stopping criterion, $\mathcal{B}$, $\mathcal{I}$, $\mathcal{K}$, $L_0$, $Q$, $T$.\\
\underline{Initialization:}
\begin{enumerate}[label={\arabic*)}]
    \item \label{alg_step: VNS step 1} Generate set $\mathcal{L} = \{\boldsymbol{p}^\ell,\, \ell = 1, \ldots, L_0\}$, where $p_i^\ell\in\mathcal{B}, \, \forall i\in \mathcal{I},\, \forall \ell = 1, \ldots, L_0$.
    \item\label{alg_step: VNS step 2} Solve $LL^k(\boldsymbol{p}^\ell), \forall \boldsymbol{p}^\ell \in \mathcal{L}, \forall k\in\mathcal{K}$. The optimal decision variables are denoted by $\boldsymbol{x}^k(\boldsymbol{p}^\ell)$.
\item \label{alg_step: VNS step 3} Evaluate the upper-level objective function $\phi^\ell:= \phi(\boldsymbol{p}^\ell, \boldsymbol{x}^k(\boldsymbol{p}^\ell))$ for all the vectors of prices $\boldsymbol{p}^\ell \in \mathcal{L}$.
\item\label{alg_step: VNS step 4} Define $\phi^{max} = \max\limits_{{\boldsymbol{p}^\ell}\in \mathcal{L}}\phi(\boldsymbol{p}^\ell, \boldsymbol{x}^k(\boldsymbol{p}^\ell))$.
\end{enumerate}
\underline{Loop:}
\begin{enumerate}[start=5,label={\arabic*)}]
    \item Set $r = 1$.
\end{enumerate}
\While{stopping criterion}

\begin{enumerate}[start=6,label={\arabic*)}]

\item Create the set $\mathcal{Q} = \left\{ \boldsymbol{p}^{\ell_q} \in\mathcal{L}, q = 1, \ldots, Q: \phi^{\ell_1} \geq \ldots \geq \phi^{\ell_Q} \right\}$.
\item Create the set $\mathcal{T} = \left\{\boldsymbol{p}\in N_r(\boldsymbol{\tilde{p}}) \text{ for }\boldsymbol{\tilde{p}}\in\mathcal{Q}\right\}$.
\item \label{alg_step: VNS step 8}Update $\mathcal{L}$, i.e., $\mathcal{L}:= \mathcal{L} \cup \mathcal{T}$.
\item Solve $LL^k(\boldsymbol{p}), \forall \boldsymbol{p} \in \mathcal{T}, \forall k\in\mathcal{K}$.
\item \label{alg_step: VNS step 10} Evaluate the upper-level objective function $\phi(\boldsymbol{p}, \boldsymbol{x}^k(\boldsymbol{p})), \forall \boldsymbol{p}\in \mathcal{T}$.
\end{enumerate}
\If{$\exists\, \boldsymbol{p}\in \mathcal{T} \text{ such that } \phi(\boldsymbol{p}, \boldsymbol{x}^k(\boldsymbol{p})) >\phi^{max}$}  $\phi^{max} := \phi(\boldsymbol{p}, \boldsymbol{x}^k(\boldsymbol{p}))$.
\Else \text{ } $r:= r+1$.
    \EndIf
\EndWhile

\vspace*{0.1cm}
\State 11) The vector of prices associated to $\phi^{max}$ is denoted by $\boldsymbol{p}^{max}$.\\

\textbf{Output: } $\left(\boldsymbol{p}^{max}, \phi^{max}\right)$.
\end{algorithmic}
\end{algorithm}

Let us show how our approach works in the previous illustrative example. Regarding the input parameters, we have already defined the sets $\mathcal{B}$, $\mathcal{I}$ and $\mathcal{K}$ in Section \ref{subs: base algorithms}. We decided to set the values $L_0=10$, $Q = 3$ and $T=6$. As a stopping criterion, we choose to run the algorithm until two iterations have been carried out. Figure \ref{fig: new points illustrative example} shows the initial set of points, $\mathcal{L}$ depicted with circles. As in the previous figure, the color of the points represents the objective value. The largest objective value so far is $\phi^{max}=228$ and is attained at point (27, 66). In the first iteration of our proposed algorithm, we create the set $\mathcal{Q}$ with three points. The points in Figure \ref{fig: new points illustrative example} that are rounded with a black circumference are those that belong to $\mathcal{Q}$. In other words, $\mathcal{Q} = \{(27, 66), (50, 18), (66, 27)\}$. In the next step, we generate the set $\mathcal{T}$ of $6$ points randomly sampled from the neighborhood of radius 1 of the previous points. Figure \ref{fig: new points illustrative example} depicts these new points in $\mathcal{T}$ with a square. Among these new sampled points, we have the true optimum, which corresponds to $(50, 34)$ with an objective value of $236$. Hence, since there exists a point with a higher objective value ($236>228$), then we now set $\phi^{max} = 236$. In practice, we do not know which point is the optimum, so we have to continue the process. We begin the second iteration with a set $\mathcal{L}$ with 16 points. Firstly, we build the set $\mathcal{Q}\subset \mathcal{L}$ that contains the three vectors of prices with the best objective values so far. Figure \ref{subfig: new points second iteration illustrative example} depicts with circles the 16 vectors of prices belonging to $\mathcal{L}$ and marked with black circles the three points of $\mathcal{Q}$. Then, we have that $\mathcal{Q} = \{(50, 34), (50, 27), (66, 27)\}$. The new 6 points randomly sampled in the neighborhood of radius 1 from the points in $\mathcal{Q}$ are depicted with squares in Figure \ref{subfig: new points second iteration illustrative example}. Since we reached the maximum number of iterations, then we finish the algorithm with the output $\phi^{max} = 236$ and $\boldsymbol{p}^{max}= (50, 34)$.

\begin{figure}[htb!]
\centering
\includegraphics[scale=0.6]{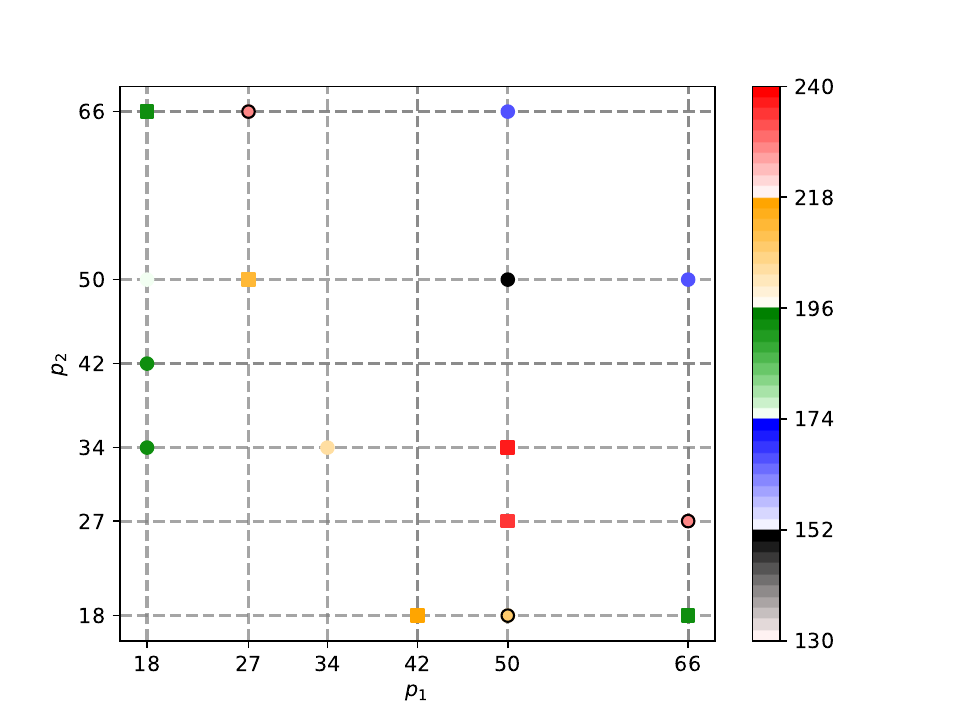}
\caption{Sets $\mathcal{L}$, $\mathcal{Q}$ and $\mathcal{T}$ for the first iteration of the illustrative example.}
 \label{fig: new points illustrative example}
    \end{figure}

\begin{figure}[htb!]
\centering
\includegraphics[scale=0.6]{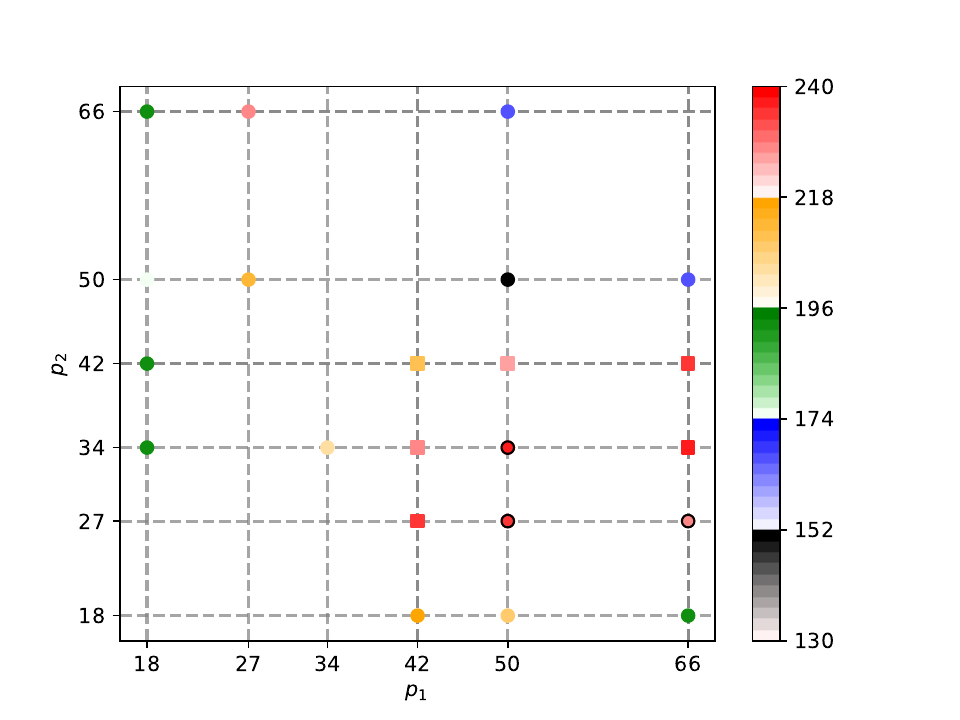}
\caption{Sets $\mathcal{L}$, $\mathcal{Q}$ and $\mathcal{T}$ for the second iteration of the illustrative example.}
 \label{subfig: new points second iteration illustrative example}
    \end{figure}

\subsubsection{The genetic-based approach}\label{subsubs: genetic approach}

We have discussed at the beginning of Section \ref{subsubs: VNS approach} that the discrete nature of the Rank Pricing Problem may cause that small changes in the vector of prices lead to big changes in the objective values. This is, in fact, why the algorithm VNS may be appropriate for this situation. Nevertheless, depending on the RPP instance we consider, its feasible region may have flat zones in terms of the objective value. This implies that exploring such a feasible region with the VNS neighborhoods may require the evaluation of many price vectors, and hence solving many optimization problems, whose solution will provide, in the best case, a slight improvement in the objective value.

To see this, consider the toy example from Figure \ref{fig: feasible solutions illustrative example}. If we focus on point $(18, 27)$ which has an objective value close to 174 (due to the light green color), we observe that four out of the five points that belong to the VNS neighborhood of radius 1 has similar objective values. Indeed, there are two of them, namely $(18, 18)$ and $(27, 18)$, with objective values less than 174. In addition, only one point out of five, $(27, 34)$, belongs to the next color range (yellow) of the incumbent point. Something similar happens if we move to the VNS neighborhood of radius 2. We have to increase the radius of the neighborhood until $r = 3$ to find a point belonging to the red range of values (which are the ones we are most interested in because these are the ones associated with  higher objective values). In other words, in this particular case, if we initialize the algorithm from (18, 27), we have to evaluate until 19 price vectors, and solving their corresponding optimization problems to obtain an objective value in the red zone. The genetic algorithm is an alternative to the VNS that allows to escape from the flat zones of the RPP feasible region without the need of exploring so many points. This does not mean that the genetic method is always better than the VNS. In fact, the numerical experiments of Section \ref{sec: numerical experience} empirically shows that this is not the case.

Now, let us explain how the basic version of the genetic algorithm works. We assume given an initial population. Then, two main steps are performed, the crossover and the mutation operators. In the crossover, two vectors are chosen from the initial population. These vectors will play the role of parents, and will be used to create the so-called child vector. Then, in the mutation phase some of the components of the child will be changed, building a new point. This process is repeated until a stopping criterion is reached.

Therefore, it is essential to define the details of how the crossover and mutation phases are done. In the context of the RPP, the genetic algorithm has been already applied in \citep{calvete2024evolutionary}. Hence, we utilize the same strategy the authors define in this paper, and just adapt it to the notation of this manuscript. Given two parents from the current population, each component of the child vector $\boldsymbol{p}^{child}$ is randomly chosen from one of the parents with equal probability. Then, in the mutation phase, each component of the price vector  $\boldsymbol{p}^{child}$ is changed with probability $1/I$ to any of the remaining possible budget values in $\mathcal{B}$.

For completeness, we describe in what follows how the genetic algorithm is applied in the context of the RPP. Moreover, in order to avoid repetitions, we will refer to the VNS pseudocode in Algorithm \ref{alg: VNS} in those parts of the genetic procedure that are very similar. We begin with an initial set $\mathcal{L}$ of $L_0$ vectors, as well as with the best objective value, $\phi^{max}$ in $\mathcal{L}$ obtained following the steps \ref{alg_step: VNS step 2} - \ref{alg_step: VNS step 4} from Algorithm \ref{alg: VNS}. Next, at each iteration, we select, from $\mathcal{L}$, the $Q$ vectors with the highest objective values, and save them in the set $\mathcal{Q}$. Now, the set $\mathcal{T}$ is to be built. It is formed by $T$ vectors that are created using the crossover and mutation operator. For each $t = 1, \ldots, T$, we randomly choose two parents from $\mathcal{Q}$, $\boldsymbol{p}^{q_1}$ and $\boldsymbol{p}^{q_2}$ to run the crossover operator and create $\boldsymbol{p}^{child}$ as follows: The component $p_i^{child}$ is selected with equal probability between $p_i^{q_1}$ and $p_i^{q_2}$ for $i\in \mathcal{I}$. Then, the mutation operator is executed. Each component $p_i^{child}$ is changed with probability $1/I$ to any possible value in $\mathcal{B}$. The resulting vector is denoted by $\boldsymbol{p}^{mut}$, and is included in the set $\mathcal{T}$. Afterward, steps \ref{alg_step: VNS step 8} - \ref{alg_step: VNS step 10} for Algorithm \ref{alg: VNS}. In the following step, we just check whether there exists a vector $\boldsymbol{p}\in \mathcal{T}$ such that $\phi(\boldsymbol{p}, \boldsymbol{x}^k(\boldsymbol{p}))$ is greater than the current maximum $\phi^{max}$. In such a case, we update $\phi^{max}$ as $\phi(\boldsymbol{p}, \boldsymbol{x}^k(\boldsymbol{p}))$. When the algorithm finishes, we just denote as $\boldsymbol{p}^{max}$ the vector of prices associated with the objective value $\phi^{max}$. Algorithm \ref{alg: genetic} shows the details of the pseudocode.

{
\renewcommand{\thealgorithm}{\textcolor{black}{\arabic{algorithm}}} 
\makeatletter
\def\ALG@name{\textcolor{black}{Algorithm}} 
\makeatother
\begin{algorithm}
\caption{\emph{genetic}}
\label{alg: genetic}
{\color{black}
\begin{algorithmic}
\State \textbf{Input:} stopping criterion, $\mathcal{B}$, $\mathcal{I}$, $\mathcal{K}$, $L_0$, $Q$, $T$.\\
\State \underline{Initialization:}
\begin{enumerate}[label={\arabic*)}]
    \item Run the steps \ref{alg_step: VNS step 1} - \ref{alg_step: VNS step 4} of Algorithm \ref{alg: VNS}.
\end{enumerate}
\underline{Loop:}
\While{stopping criterion}
\begin{enumerate}[label={\arabic*)}, start = 2]
\item Create the set $\mathcal{Q} = \left\{ \boldsymbol{p}^{\ell_q} \in\mathcal{L}, q = 1, \ldots, Q: \phi^{\ell_1} \geq \ldots \geq \phi^{\ell_Q} \right\}$.
\item $\mathcal{T} = \emptyset$.
\end{enumerate}
\For{$t \in \{1, \ldots, T\}$}\\
 \begin{enumerate}[label={\arabic*)}, start = 4, leftmargin=1.5cm]
 \item Select $\boldsymbol{p}^{q_1}$ and $\boldsymbol{p}^{q_2}$ from $\mathcal{Q}$.
 \item Create $\boldsymbol{p}^{child}$ from $\boldsymbol{p}^{q_1}$ and $\boldsymbol{p}^{q_2}$ in this way:
 \end{enumerate}
\begin{equation*}
    p_i^{child} =\beta*p_i^{q_1}+ (1-\beta)*p_i^{q_2},\, i = 1, \ldots, I
\end{equation*}
\hspace*{1.5cm} where $\beta$ follows a Bernouilli distribution with probability $0.5$.
\begin{enumerate}[label={\arabic*)}, start = 6, leftmargin=1.5cm]
\item Create the point $\boldsymbol{p}^{mut}$ from $\boldsymbol{p}^{child}$, where $p^{mut}_i$ changes with probability $1/I$, $\forall i$ to any possible budget in $\mathcal{B}$. 
\item Update $\mathcal{T}$, i.e., $\mathcal{T}:= \mathcal{T} \cup \{\boldsymbol{p}^{mut}\}$.
\end{enumerate}
\EndFor
\begin{enumerate}[label={\arabic*)}, start = 8, , leftmargin=0.9cm]
\item \label{alg_step: gen step 8} Run steps \ref{alg_step: VNS step 8} - \ref{alg_step: VNS step 10} of Algorithm \ref{alg: VNS}.
\end{enumerate}

\If{$\exists\, \boldsymbol{p}\in \mathcal{T} \text{ such that } \phi(\boldsymbol{p}, \boldsymbol{x}^k(\boldsymbol{p})) >\phi^{max}$} $\phi^{max} := \phi(\boldsymbol{p}, \boldsymbol{x}^k(\boldsymbol{p}))$.
    \EndIf
\EndWhile
\begin{enumerate}[label={\arabic*)}, start = 9]
\item The vector of prices associated to $\phi^{max}$ is denoted by $\boldsymbol{p}^{max}$.
\end{enumerate}

\State \textbf{Output: } $\left(\boldsymbol{p}^{max}, \phi^{max}\right)$.
\end{algorithmic}
}
\end{algorithm}
}

\subsection{Initialization procedure} \label{subs: initialization procedure}

One of the assumptions of the heuristic algorithms described in Section \ref{subs: base algorithms} is the availability of an initial set $\mathcal{L}$ consisting of $L_0$ price vectors. We have considered two approaches for generating this set. The first relies on a standard random selection. More specifically, for each product $i = 1, \ldots, I$, its price, $p_i$ is randomly chosen from the set of sorted budgets, $\mathcal{B}$, following a discrete uniform probability distribution. With this process, we obtain a vector $\boldsymbol{p}$ that is to be included in the initial set $\mathcal{L}$. This procedure is repeated $L_0$ times. The second approach, proposed in \citep{calvete2024evolutionary}, constructs the first vector greedily: customers are sorted in decreasing order based on their budget, and for each customer, their most preferred available product is priced at their budget value. Our numerical experiments indicate that this initial greedy vector often leads to strong objective values, particularly when the number of products is comparable to the number of customers. The subsequent vectors in the population are randomly sampled from various discrete distributions.

\subsection{Local improvements} \label{subs: local improvements}

One of the key aspects of \citep{calvete2024evolutionary} for solving the Rank Pricing Problem is the local search method applied immediately after executing the genetic algorithm. Although their approach yields good results, it requires solving a potentially large number of optimization problems. In contrast, this paper introduces four new local search enhancements that leverage the problem's structure, thereby avoiding the need to solve additional optimization problems. Section \ref{subsubs: optimization-based local search} outlines the main features of the local search proposed in \citep{calvete2024evolutionary}, while Sections \ref{subsubs: slack local search}–\ref{subsubs: binary local search} detail our new approaches. The core idea behind these improvements is to adjust an incumbent price vector to obtain another that enhances the upper-level objective function. These adjustments follow different criteria, which are explained below.

\subsubsection{Optimization-based local search} \label{subsubs: optimization-based local search}
This section explains the local improvement proposed by \citep{calvete2024evolutionary}. Starting for one of the vectors $\boldsymbol{p}^{mut}$ obtained after running the genetic approach, the products are randomly sorted. Then, for each product, they modify the current price to one of the possible available budgets, and compute the new revenue associated to the vector, by solving the lower-level problems \eqref{eq: lower level} and evaluating the upper-level objective function \eqref{eq: upper level objective function}. If the new objective value is better, then the new vector of prices is kept. Otherwise, the procedure continues with the next budget value and the next product. This process is repeated for all the points of the population obtained with the genetic algorithm. Let us show with the illustrative example of Section \ref{sec: methodology} how the optimization-based local search works. Imagine that after running the genetic algorithm we have the point $\boldsymbol{p}^{mut} = (42, 34)$ with an objective value of 228 (coming from $3\cdot 42 + 3 \cdot 34$). Then, assume that we sort the products in the order $2 -1$. In other words, we first change product $2$, and next, product $1$ is modified. In the following step, we assign to product $2$ a different price among those available. For instance, we choose to set $p_2$ to $27$, yielding the vector $(42, 27)$. After solving $LL^k((42, 27)),\, \forall k$, and evaluating the upper-level objective function we get a value of 234. Since $234>228$, we keep the new vector, and continue the process.

\subsubsection{Slack local search} \label{subsubs: slack local search}

The four local searches that we propose in this paper do not need to solve new optimization problems but just take advantage of the structure of the problem.

The strategy proposed in this section keep intact the associated binary variables. In other words, the binary decision variables associated to the vector before performing the local search are exactly the same as the ones of the vector of prices yielded after the local improvement. The slack local search works as follows: given a vector of prices $\boldsymbol{p}$, if there exists a product $i$ that is chosen by a set of customers $\hat{\mathcal{K}}\subset \mathcal{K}$ such that the budgets of the customers in $\hat{\mathcal{K}}$ are always larger than $p_i$, i.e., $b^k>p_i, \forall k\in \hat{\mathcal{K}}$, then we can change the price of product $i$, $p_i$ to the minimum budget value of the customer in $\hat{\mathcal{K}}$. In other words, we update $p_i:=\min\limits_{k \in \hat{\mathcal{K}}} b^k$. This change implies an improvement in the objective value and no changes in the binaries because all the customers will choose exactly the same products as before, but just at a different price. Let us focus on Table \ref{tab: data illustrative example} to illustrate how it works. Let us take point $(34, 34)$. Product 1 has a price of 34 and is chosen by customers 2, 6 and 7. Their budget values are $b^2 = 66$, $b^6= 50$ and $b^7 = 42$. On the other hand, product 2 has also a price of 34, and is selected by customers 4, 5 and 8, with budgets 34, 66 and 42, respectively. The upper-level objective function is $3\cdot34 + 3\cdot 34 = 204$. We observe that product $1$ has some \emph{slack} in terms of the proposed price. Indeed, we have that $b^k>34, \forall k \in \hat{\mathcal{K}}$, where $\hat{\mathcal{K}} = \{2, 6, 7\}$.  The customers in $\hat{\mathcal{K}}$ can still choose product 1 if the price is 42, i.e., the minimum of $b^k, k \in \hat{\mathcal{K}}$, yielding point (42, 34). This change has no modifications in the binary variables because all the customers get the same products as before, but has a modification in the objective function. In particular, it goes from 204 to the value $3\cdot 42 + 3\cdot 34 = 228$.

\subsubsection{Fill local search}

This is the first proposed strategy that imply changes in the binary variables. Assume given a vector of prices $\boldsymbol{p}$, and assume that there exists a product $i^*$ that has not been chosen by any customer because the price $p_{i^*}$ is too high. Equivalently, $\exists i^* \in \mathcal{I}: x_{i^*}^k = 0, \forall k\in\mathcal{K}$. The proposed strategy to \emph{fill} those gaps works as follows: let us denote by $\hat{\mathcal{K}}$ the set of customers that do not choose any product, that is to say, for a given $k\in\hat{\mathcal{K}}$ we have that $x_i^k= 0, \forall i$. Now, we set to the price $p_{i^*}$ the minimum budget value of the customers in $\hat{\mathcal{K}}$, i.e., $p_{i^*}:=\min\limits_{k\in \hat{\mathcal{K}}} b^k$. This way, the binary variables changed as follows: $x_{i^*}^k = 1, \forall k \in \hat{\mathcal{K}}$, whereas $x_{i^*}^k$ still takes the value zero for the customers $k\in\mathcal{K}\setminus \hat{\mathcal{K}}$. Besides, the objective value has been increased because the customers in $\hat{\mathcal{K}}$ that did not contribute to the objective function, now they sum the new assigned price. Let us go back to the Table \ref{tab: data illustrative example} to illustrate this. Before starting the example two minor changes should be done. Assume that $s_1^5 = 2$ and $s_2^5 = None$. This means that now, customer 5 can only adquire product 1. It is easy to see that if we take the point $(66, 66)$, then product 1 is selected by customers $2$ and $5$, but product 2 is bought by nobody as the remaining customers, that is to say, the customers in $\hat{\mathcal{K}} = \{1, 3, 4, 6, 8\}$, have a budget value less than 66 which is the price fixed for this product. Therefore, $x_2^k = 0, \forall k \in \mathcal{K}$. In particular, $x_2^k = 0, \forall k \in \hat{\mathcal{K}}$. The upper-level objective function associated to the point $(66, 66)$ is $2 \cdot 66 = 132$. On the other hand, the minimum budget value for the customers in $\hat{\mathcal{K}}$ is 18. This way, since all the customers in $\hat{\mathcal{K}}$ can buy product 2 (because $s_2^k>0, \forall k \in \hat{\mathcal{K}}$), we are sure that if we set $p_2=18$, then all the customers in $\hat{\mathcal{K}}$ will take it. The binary variables will be $x_2^k = 1$ for $k\in \hat{\mathcal{K}}$ and $x_2^k= 0$ for $k\in \mathcal{K}\setminus \hat{\mathcal{K}}$, and the new objective function increases till the value $2\cdot 66 + 6 \cdot 18 = 240$.

\subsubsection{Reassigment local search}\label{subsubs: TSP local search}

The local improvement performed with this strategy also implies a modification in the binary decision variables.  The approach is done as follows: we assume given a vector $\boldsymbol{p}$ whose prices are adjusted to the minimum budgets of the clients that acquire it. In other words, we have, for a product $i$ that $p_i= \min\limits_{
k \in \mathcal{K}^i} b^k$, where $\mathcal{K}^i\subset \mathcal{K}$ is the set of customers that select to buy product $i$. Note that, due to constraint \eqref{eq: lower level first constraint}, we have that $\mathcal{K}^{i_1} \cup \mathcal{K}^{i_2} = \emptyset$ for $i_1\neq i_2$.  We then ask these questions: What happens if we \emph{reassign} the price $p_i$ to the second cheapest budget (instead of the cheapest one, as it is done so far)? Does the objective function increase its value? Which product will be assigned to the customer with the cheapest budget? Let us check it with the illustrative example of Table \ref{tab: data illustrative example}. Consider the point $(18, 27)$. Note that customers $1, 2, 6$ and $7$ with budgets $18, 66, 50$ and 42, respectively, will select product 1. Equivalently, $x_1^1 = x_1^2 = x_1^6 = x_1^7 = 1$, whereas $x_1^3 = x_1^4 = x_1^5= x_1^8 =0$. On the other hand, product 2 with price $p_2 = 27$ is chosen by customers $3, 4, 5$ and $8$. The objective function has a value of $4 \cdot 18 + 4\cdot 27 = 180$. What does it happen if we change the price of product 1 to the budget of the second poorest customer? That is to say, if we change $p_1 = 18$ to $p_1 = 42$? In this particular case, customer 1 takes no product, i.e., $x_1^1 = x_2^1 = 0$, thus the binary variables change. In contrast, the remaining customers mantain their assignation. The objective function increases from 180 to the value of $3\cdot 42 + 4\cdot 27 = 234$. Note that, contrary to what happens in the previous local improvement strategies, this procedure may not imply a positive result. In other words, it may occur that changing the price of a product to the budget of the second poorest customer decreases the objective function. If this is the case, then, we keep the original point. Nevertheless, the computational experience of Section \ref{sec: numerical experience} shows that, in general, this approach yield good results.

\subsubsection{Conditional reassignment local search} \label{subsubs: binary local search}

This approach is similar to the one presented in Section \ref{subsubs: TSP local search}. The main difference is that the local changes done with this strategy always yield better objective values, contrary to what happened before, where worse scenarios are possible. 

We assume given a vector of prices $\boldsymbol{p}$, where the prices have no slack, that is to say, the minimum budget of the customers that choose each product $i$ coincides with the price assigned to product $i$, $p_i$. Therefore, given a product $i$, and a customer $k^*$ such that $p_i = b^{k^*}$, the idea is to find if there exists another product $\hat{i}$ that can be selected by customer $k^*$ ($s_{\hat{i}}^{k^*}>0$) and whose price coincides with that of product $i$, i.e., $p_{\hat{i}} = b^{k^*}$. If this is the case, we set the price of product $i$ to the second cheapest budget. This way, all the customers, except the poorest one, $k^*$, that originally are assigned to product $i$ still select this product, whereas the poorest customer, $k^*$ is assigned to product $\hat{i}$ whose price coincides with the original one. This change in the prices causes modifications in the binary variables. Let us show how it works with the illustrative example of Table 1. Assume given the point $(42, 42)$. We can see that product 1 is assigned to the customers 2, 6 and 7 with budgets, 66, 50 and 42, respectively. Hence, $x_1^2 = x_1^6 = x_1^7 = 1$ and $x_1^1 = x_1^3 = x_1^4= x_1^5 = x_1^8 =0$. On the other hand, product 2 is assigned to customers $5$ and $8$, and consequently, $x_2^5 = x_2^8 = 1$ and $x_2^1 = x_2^2 = x_2^3 = x_2^4= x_2^6 = x_2^7 =0$. The objective function takes the value $3\cdot 42 + 2 \cdot 42 = 210$. It can be observed that if price of product $1$ is assigned to the second cheapest budget, i.e., $p_1 = 50$, then customer 7 (who was the one with budget 42) now selects product 2, yielding $x_1^7 = 0$ and $x_2^7 = 1$. The assignation of the rest of customers remain unchanged. Finally, the objective function associated to the new point $(50, 42)$ is $2\cdot 50 + 3\cdot 42 = 226$, thus increasing the objective value.

To sum up this section, we include in Figure \ref{fig: flowchart} a flowchart of the different heuristic methodologies explained in this paper. The first dashed box describes the initial step, whereas the second dashed box explains how the solutions are updated in the loop.

\begin{figure}
    \centering
\begin{tikzpicture}
    \node (generateL) [process] {Generate $\mathcal{L}$ (random/greedy)};
    \node (solveLL) [process, below of=generateL, yshift=-1.5cm] {Solve $LL^k(\cdot)$ and evaluate the upper-level function $\phi(\cdot, \cdot)$};    
    
    \node (createQ) [process, below of=solveLL, yshift=-3cm,align=center] {Create $\mathcal{Q}$ \\ Create $\mathcal{T}$ (VNS/genetic)};
    
    \node (updateL) [process, below of=createQ, yshift=-1.5cm] {Update $\mathcal{L}$};
    
    \node (solveIter) [process, below of=updateL, yshift=-1.5cm] {Solve $LL^k(\cdot)$ and evaluate the upper-level function $\phi(\cdot, \cdot)$};    
    
    \node (localSearch) [process, below of=solveIter, yshift=-1.5cm, align=center] {Local search\\(opt-based/slack/fill/\\reassignment/cond. reassignment)};
        
    \draw [arrow] (generateL) -- (solveLL);
    
    \draw [arrow] (solveLL) -- (createQ);
    
    \draw [arrow] (createQ) -- (updateL);
    \draw [arrow] (updateL) -- (solveIter);
    \draw [arrow] (solveIter) -- (localSearch);   
    \draw [arrow] (localSearch.south) |- ++(-2.5,-0.3) -| ++(-3,4.5) |- (createQ.west);    
    
    \node[draw, dashed, fit=(generateL)(solveLL), inner sep=0.75cm] (init) {};
    
    \node at (init.north west) [above] {Initialization};
    
    \node[draw, dashed, fit=(createQ)(localSearch)(solveIter), inner sep=0.75cm] (loop) {};
    
    \node at (loop.north west) [above] {Loop};
    
\end{tikzpicture}
\captionsetup{labelfont={color=black}}
\caption{Flowchart of the different heuristic methodologies.}
\label{fig: flowchart}
\end{figure}

}
\section{Numerical Experiments} \label{sec: numerical experience}

\subsection{Experimental Setup}
In this section, we explain the different computational experiments carried out in this paper. We test our approach in three instances, where we modify the number of customers, $K$, and products, $I$. In particular, we apply our proposal on the instances with $(K, I)$ belonging to the set $\{(30, 5), (30, 25), (60, 50)\}$. The maximum size of the instances solved in \citep{dominguez2021rank} is $K = 30$ and $I = 25$. We decided to include a new instance with a larger size,  $(K, I) = (60, 50)$, to test how our proposal works in this example. The data for the three instances can be downloaded from \citep{OASYS2024RPP}. 

Let us notate the different strategies compared in this section. Firstly, we denote by \emph{exact} and \emph{naive} the approaches detailed in Sections \ref{subs: exact approach} and \ref{subs: naive}, respectively. Secondly, we compare the proposal with the strategy given in \cite{calvete2024evolutionary} that we denote as \emph{benchmark}. This methodology is obtained after running the genetic algorithm of Section \ref{subsubs: genetic approach} with the greedy initialization process the authors propose in that paper and improved with the optimization-based local search of Section \ref{subsubs: optimization-based local search}. The optimization-based local improvement is run after step \ref{alg_step: gen step 8} in Algorithm \ref{alg: genetic}.

The proposed approaches where a heuristic, VNS or genetic, is improved with four local searches are denoted, respectively, \emph{VNS + greedy + sfrc} and \emph{genetic + greedy + sfrc}. In the case of \emph{VNS + greedy + sfrc}, it means that we run the VNS algorithm of Section \ref{subsubs: VNS approach} using the greedy initialization and applying, in this order, the local searches \emph{slack}, \emph{fill}, \emph{reassignment} and \emph{conditional reassignment} from Sections \ref{subsubs: slack local search} - \ref{subsubs: binary local search}. Note that these local searches are run sequentially one after another. In addition, when the VNS applies the local improvements are executed at each iteration after evaluating the upper-level objective function, that is to say just after running step \ref{alg_step: VNS step 10} of Algorithm \ref{alg: VNS}. The vector of prices and the associated objective values that are computed in the following steps correspond to those obtained with the local searches. On the other hand, if the \emph{gen + greedy + sfrc} approach is executed, then we run the genetic algorithm from Section \ref{subsubs: genetic approach} initialized with the greedy strategy and improved with the four proposed local searches of Sections \ref{subsubs: slack local search} - \ref{subsubs: binary local search}. In this case, the four local searches are run after executing Step \ref{alg_step: gen step 8} of Algorithm \ref{alg: genetic}.

In order to study the effect of the greedy initialization and the four local searches in our proposal, we also analyse four different approaches, denoted as \emph{VNS}, \emph{genetic}, \emph{VNS + greedy} and \emph{gen + greedy}. If \emph{VNS} or \emph{genetic} are applied, then Algorithms \ref{alg: VNS} or \ref{alg: genetic}, are respectively run. In this case, both heuristics are initialized randomly as explained in Section \ref{subs: initialization procedure}. Nevertheless, when the methods \emph{VNS + greedy} or \emph{gen + greedy} are executed, Algorithms \ref{alg: VNS} or \ref{alg: genetic} are respectively run, but initialized with the greedy approach.

To avoid confusion in the different strategies applied, we include Table \ref{tab: summary heuristics}. Columns 2, 3 and 4 of this table indicate, respectively, the initialization procedure, the base algorithm and the local improvement used in the heuristic methodologies \emph{benchmark}, \emph{VNS}, \emph{VNS + greedy}, \emph{VNS + greedy + sfrc}, \emph{genetic}, \emph{gen + greedy} and \emph{gen + greedy + sfrc}.  Note that the \emph{exact} and \emph{naive} approaches have not been included in the table since we found them simpler to follow.

\begin{table}[h]
    \centering
    \arrayrulecolor{black} 
    \begin{tabular}{>{\color{black}}c|>{\color{black}}c|>{\color{black}}c|>{\color{black}}c}
         & initialization &base algorithm & local improvement \\
        \hline
        \emph{benchmark}&greedy & genetic& optimization-based\\
        \hline
        \emph{VNS}& random &VNS & - \\
        \hline
        \emph{VNS + greedy}& greedy &VNS & - \\
        \hline
        \multirow{4}{*}{\emph{VNS + greedy + sfrc}} &\multirow{4}{*}{greedy} &\multirow{4}{*}{VNS} & slack + fill + \\
        && &reassignment + \\
        &&&conditional \\
        &&&reassignment\\
        \hline
        \emph{genetic}& random &genetic & - \\
        \hline
        \emph{gen + greedy}& greedy &genetic & - \\
        \hline
        \multirow{4}{*}{\emph{gen + greedy + sfrc}} &\multirow{4}{*}{greedy} &\multirow{4}{*}{genetic} & slack + fill + \\
        && &reassignment + \\
        &&&conditional \\
        &&&reassignment\\
    \end{tabular}    \captionsetup{labelfont={color=black}, textfont={color=black}}%
    \caption{Summary of the heuristic approaches.}
    \label{tab: summary heuristics}
\end{table}

For the \emph{exact} approach, we solve Problem \eqref{eq: single-level reformulation} using the off-the-shelf solver Gurobi, \citep{gurobi}. We set all the Gurobi parameters to their default value except the number of threads used. This parameter, which in Gurobi is called \texttt{Threads}, is fixed to 1. This choice was made in order not to take advantage of the different threads of the computer when running the experiments. The full list of the Gurobi hyperparameters can be found  here\footnote{\href{https://docs.gurobi.com/projects/optimizer/en/current/reference/parameters.html}{https://docs.gurobi.com/projects/optimizer/en/current/reference/parameters.html}}. The single-level reformulation \eqref{eq: single-level reformulation} is solved with a time limit of 600s. We consider the optimal solution the one given after this limit, even though the solver has not certified the optimality.  Regarding the remaining approaches, the lower-level problems of type \eqref{eq: lower level} have also been solved using Gurobi with the same settings. For the sake of comparison, these problems have been sequentially run, even though this task can be parallelized if we want to reduce the computational time even more.

We compare the approaches in terms of the objective values. Since the results for all the approaches, except the \emph{exact} may suffer variations depending on the random generation of the vector of prices, we decided to run the methodologies 1000 times to get fair conclusions. This way, a distribution of the results is shown. By default, we set the input parameters $L_0,\, Q$ and $T$ of the Algorithm \ref{alg: VNS} to the values $1000, 100$, and $500$, respectively. On the other hand, we set $L_0 = 1000$, $Q = 1000$ and $T = 500$ in the genetic strategy of Algorithm \ref{alg: genetic}.

All the experiments have been run on a Linux-based server with CPUs clocking at 2.6 GHz, 1 thread, and 8 GB of RAM using Python 3.11.4, \citep{python}. The optimization problems have been solved using Gurobi 10.0.3, \citep{gurobi}.

\subsection{Case study with 30 customers and 5 products}

This first instance is quite simple to solve in Gurobi. Indeed, even though we set a maximum time limit of 600 seconds, Gurobi is able to solve this instance in approximately 10 seconds, yielding an optimal objective value of 807. Hence, in this particular case, it makes no sense to utilize our proposal because Gurobi already provides the optimal solution in a short period of time. However, we found it interesting to compare the approach presented here with the \emph{naive} method of Section \ref{subs: naive} because, as we will see later, despite this instance being easy to solve, the \emph{naive} heuristic yields bad results compared to the proposed methodologies.

\begin{figure}[htb!]
\centering
\includegraphics[scale = 0.6]{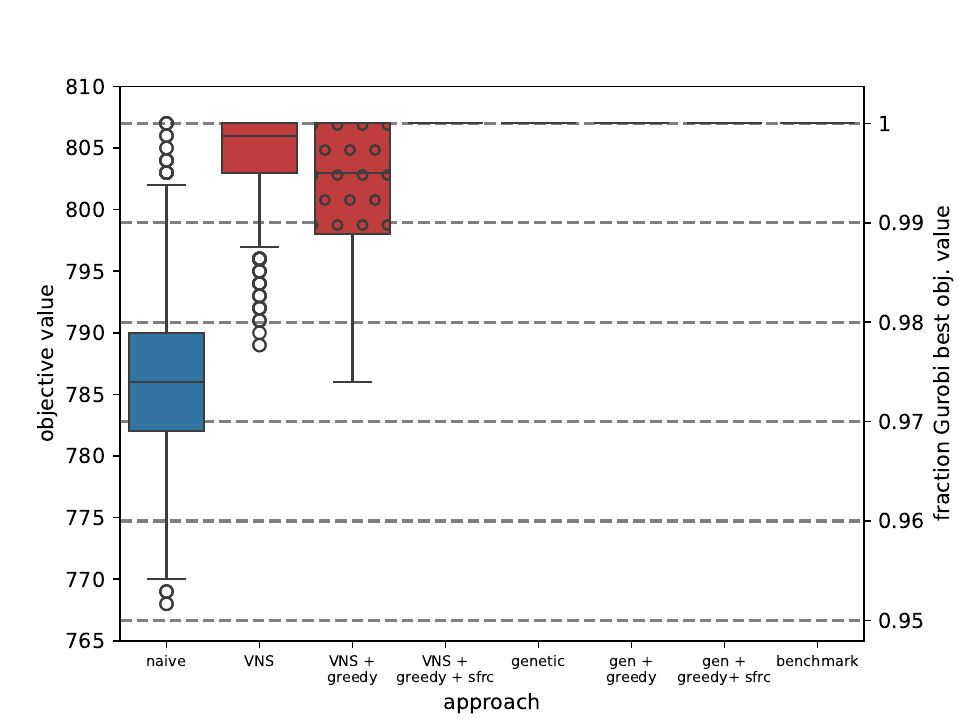}
\captionsetup{labelfont={color=black}}
\caption{Comparison of the different approaches for the instance with 30 customers and 5 products.}
 \label{subfig: boxplot max values 30c_5p}
    \end{figure}

Figure \ref{subfig: boxplot max values 30c_5p} shows the results obtained after running 1000 simulations of the approaches \emph{naive}, \emph{VNS}, \emph{VNS + greedy}, \emph{VNS + greedy + sfrc}, \emph{genetic}, \emph{gen + greedy}, \emph{gen + greedy + sfrc} and \emph{benchmark}. Note that the \emph{naive} methodology is colored in blue, whereas all the VNS heuristics are plotted in red. In addition, the genetic approaches are colored in green, and the benchmark method in orange. However, the last two colors cannot be appreciated because of the inexistent variance of the results. In addition, the output of the \emph{VNS + greedy} approach and also \emph{gen + greedy} are highlighted with circles in the boxplots. However, in the last case, it is impossible to see again due to the zero variance of the results. 

In all the cases, we use the number of points created in the set $\mathcal{L}$ as the stopping criterion. We decided to generate $|\mathcal{L}|= 24000$ points in this case. Note that the total number of feasible solutions is equal to $M^I$, where $I$ is the number of products, and $M$ is the number of unique budgets. Hence, in our case, we have $27^{5}$ solutions, which is approximately $ 10^{7}$ points to evaluate, and therefore, evaluating $24000$ points is a small quantity if we compare it with the total amount of solutions.

The left y-axis of Figure \ref{subfig: boxplot max values 30c_5p} indicates the value $\phi^{max}$ obtained with all the algorithms. On the other hand, the right y-axis shows the fraction of the optimal objective value we obtain with Gurobi after 10 seconds. For example, in this particular case, we know the optimal objective value is attained at 807. Hence, the horizontal line at 0.99 indicates the value $807\times 0.99 = 798.93$. Consequently, if the values of the boxplots are situated between 0.99 and 1, it means that those objective values are between 798.93 and 807. This is exactly what happens with our proposal, especially in the approaches \emph{VNS + greedy + sfrc, genetic, gen + greedy} and \emph{gen + greedy + sfrc}. In these cases, the heuristic methodologies are able to reach the true optimum in the 1000 simulations. The same happens in the \emph{benchmark} approach. Note however, that worse results in terms of the objective values are obtained with the approaches \emph{VNS} and \emph{VNS + greedy}, even though they are much better than the \emph{naive} strategy.
In addition, we want to compare \emph{VNS} and \emph{VNS + greedy} to analyze the effect of the different initialization methods that we use. Note that, in this instance, where the number of customers is much larger than the number of products (30 versus 5), the greedy initialization makes no benefits in terms of the objective values. In fact, the results of the approach \emph{VNS + greedy} are worse than those of \emph{VNS}, and has more variability.

In any case, as mentioned at the beginning of this section, this instance can be easily solved with Gurobi in less than 600 seconds. Then, our approach here is not useful. The following sections provide further results in more challenging instances.

\subsection{Case study with 30 customers and 25 products}

In the second case study, we analyze an instance with $K=30$ customers and $I=25$ products. As before, the time limit when solving the single-level reformulation is 600 seconds, and the number of runs of all the algorithms is 1000. We set the number of vector prices in $\mathcal{L}$ to 24000 as stopping criteria. Note that, in this instance, the number of evaluated points is even smaller than in the previous case when compared with the total number of feasible solutions, that is, $23^{25}\approx 10^{34}$.

\begin{figure}[htb!]
\centering
\includegraphics[scale = 0.7]{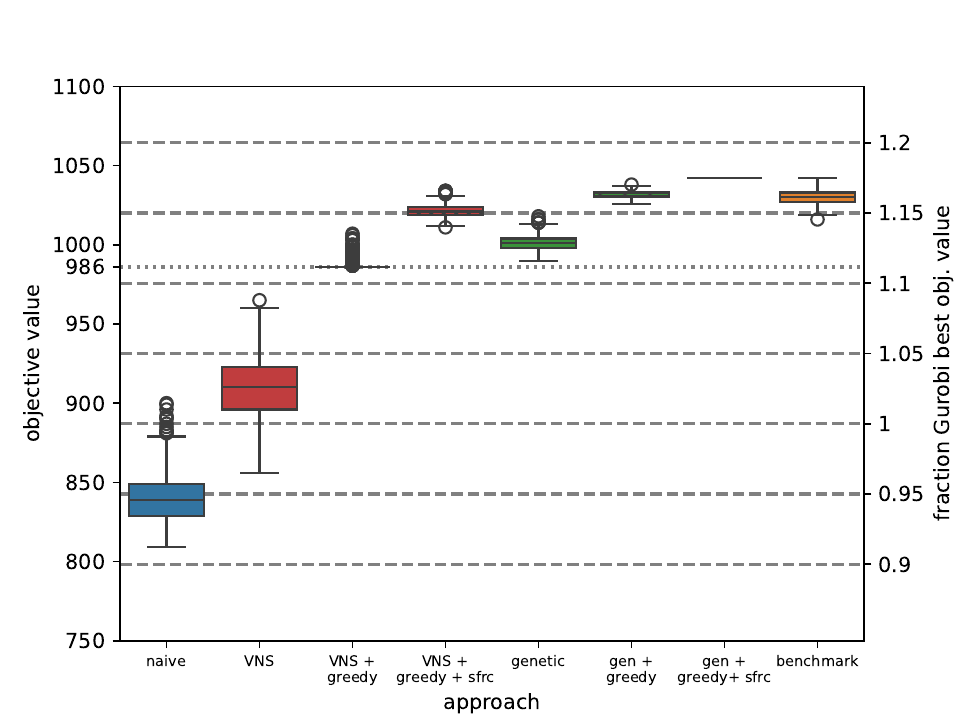}
\captionsetup{labelfont={color=black}}
\caption{Comparison of the different approaches in the instance with 30 customers and 25 products.}
\label{subfig: boxplot 30c and 25p}
\end{figure}

Figure \ref{subfig: boxplot 30c and 25p} shows the boxplots of the results obtained with the \emph{naive}, \emph{VNS}-based, \emph{genetic}-based and \emph{benchmark} approaches depicted in blue, red, green and orange, respectively. The optimal objective value provided by Gurobi in 600 seconds is 887. The dotted line at 986 represents the objective value of the first vector of prices evaluated in the greedy initialization (Section \ref{subs: initialization procedure}).

\begin{figure}[htb!]
\centering
\includegraphics[scale = 0.6]{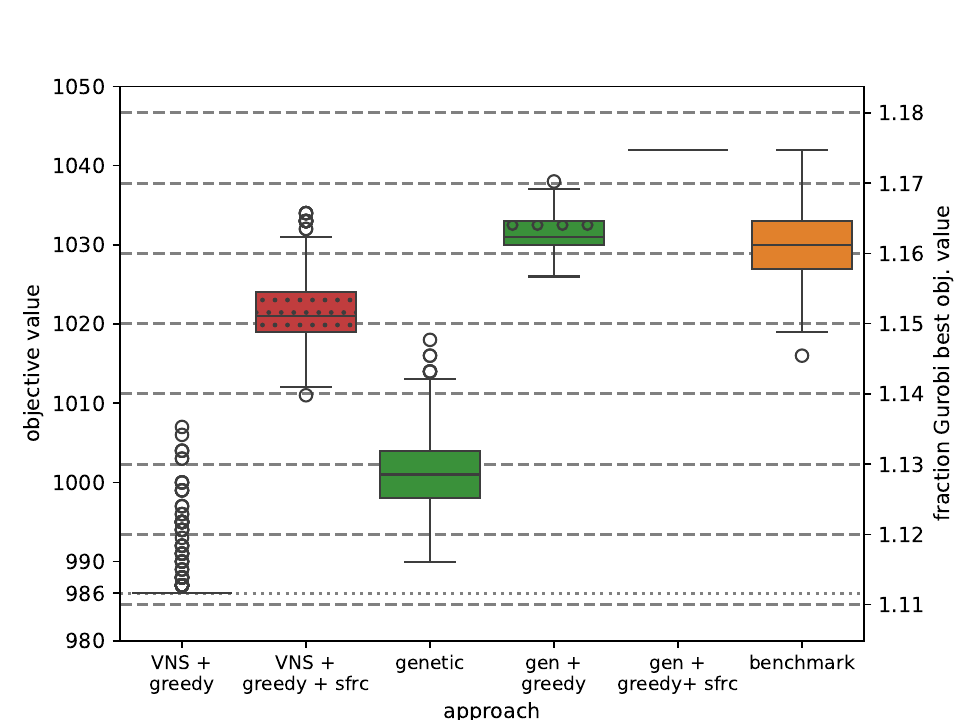}
\captionsetup{labelfont={color=black}}
\caption{Zoom of the comparison of the different approaches in the instance with 30 customers and 25 products.}
\label{fig: boxplot 30c and 25p zoom}
\end{figure}

We can observe that the outliers of the boxplot for the \emph{naive} method correspond to the cases where it reaches the best integer solution found by Gurobi. In other words, the \emph{naive} method rarely reaches this optimum. Indeed, the (central) $50\%$ of the runs of the \emph{naive} approach yield an objective value close to 840, that is, close to 0.95 times the optimum. On the contrary, the \emph{VNS} proposal got larger results than those provided by Gurobi in more than 75\% of the runs. In addition, the remaining proposed strategies, namely \emph{VNS + greedy}, \emph{VNS + greedy + sfrc}, \emph{genetic}, \emph{gen + greedy} and \emph {gen + greedy + sfrc} yields results which are larger than 1.1 times the optimal value obtained by Gurobi. This effect also occurs in the \emph{benchmark} approach. Note, however, that, in this instance, the best outputs are obtained with the approach \emph{gen + greedy + sftb}. Particularly, the 1000 simulations always attain the same value, which is 1042, and that we have empirically tested that is the true optimum of this example. For a better comparison of the results, we plot in Figure \ref{fig: boxplot 30c and 25p zoom} those boxplots with objective values larger than 1.1 times the optimal value obtained by Gurobi, i.e., larger than $1.1 \times 887 = 975.5$.

Besides, regarding the initialization procedure, we observe that, in this instance, where the number of products is similar to the number of customers (30 versus 25), the greedy initialization has a positive effect. This is particularly evident if we compare the approaches \emph{VNS} and \emph{VNS + greedy}. In addition, the comparison of the approaches without and with the four local searches \emph{slack}, \emph{fill}, \emph{reassignment} and \emph{conditional reassignment} are also noteworthy. Whether we focus at the \emph{VNS} or \emph{genetic} counterparts, we observe that the objective values obtained with \emph{VNS + greedy + sfrc} or \emph{gen + greedy + sfrc} are better than those of \emph{VNS + greedy} or \emph{gen + greedy}, respectively. Moreover, if we compare the approach \emph{gen + greedy + sfrc}, which uses our proposed local searches, with the \emph{benchmark} strategy, that applies an optimization-based local improvement, we see that the \emph{benchmark} methodology only attains the true optimum, 1042, in some cases, whereas our proposal always reaches this point.

\begin{figure}[htb!]
\centering
\includegraphics[scale = 0.6]{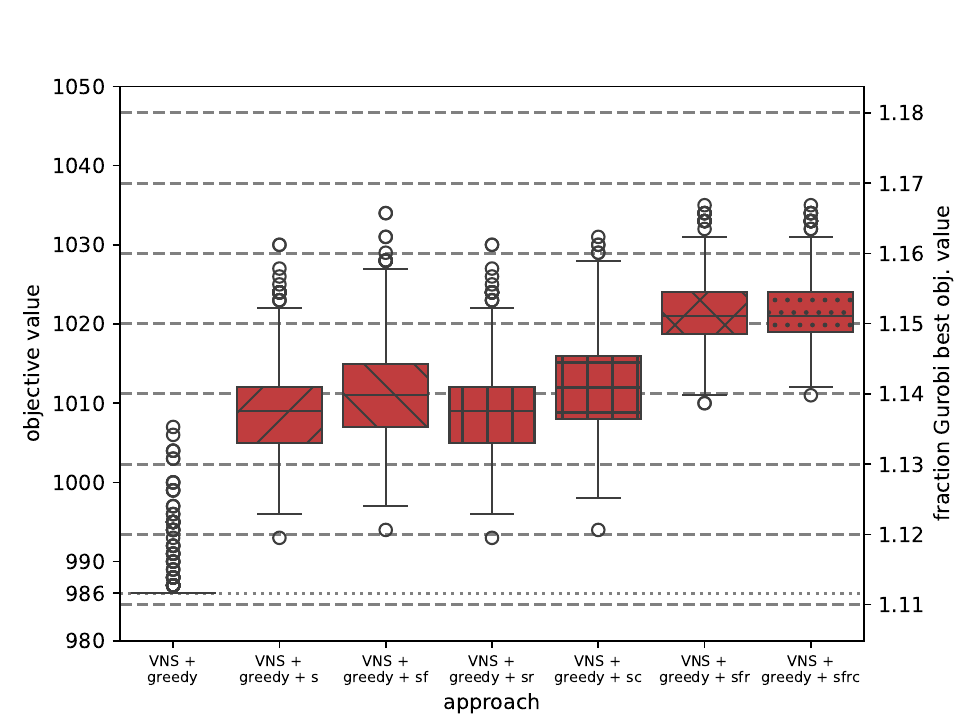}
\caption{Comparison of the four local searches using the VNS approach in the instance with 30 customers and 25 products.}
\label{fig: boxplot 30c and 25p sensitivity local search VNS zoom}
\end{figure}

\begin{figure}[htb!]
\centering
\includegraphics[scale = 0.6]{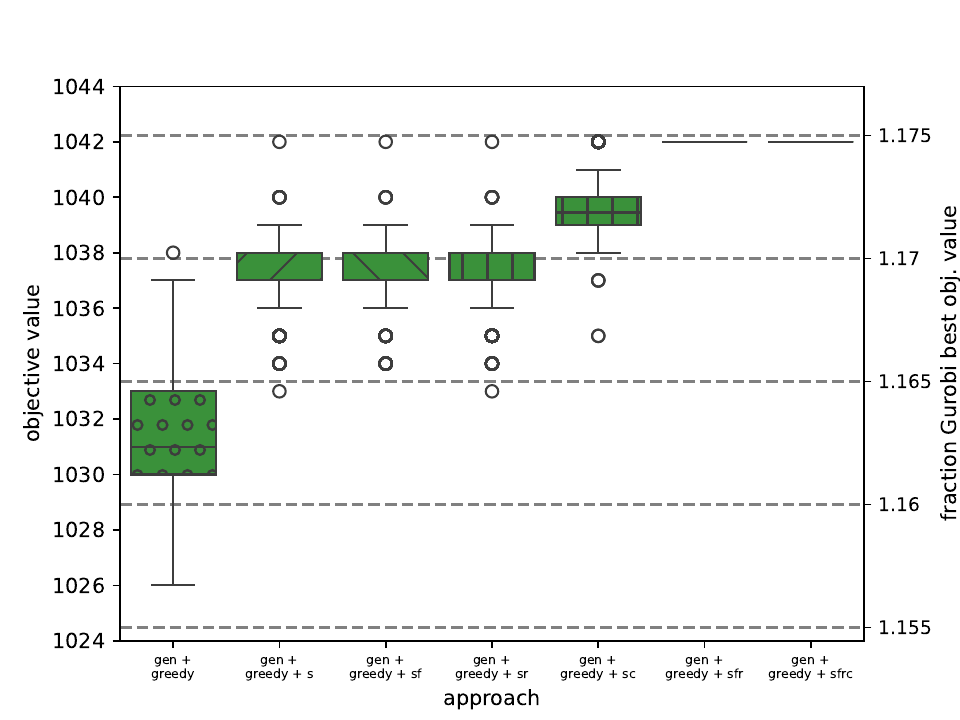}
\caption{Comparison of the four local searches using the genetic approach in the instance with 30 customers and 25 products.}
\label{fig: boxplot 30c and 25p sensitivity local search gen zoom}
\end{figure}

In what follows, let us analyze the sensitivity of the proposed local searches for the VNS-based and genetic-based heuristics. Figure \ref{fig: boxplot 30c and 25p sensitivity local search VNS zoom} shows the results in the VNS case, whereas Figure \ref{fig: boxplot 30c and 25p sensitivity local search gen zoom} depicts the analogous results for the genetic approach. In these experiments, we run 1000 times the default base heuristic with the greedy initialization, i.e. \emph{VNS + greedy} or \emph{gen + greedy} and different variants where the local searches are compared. For instance, in the case of the VNS algorithm, if we run \emph{VNS + greedy + s}, it means that the algorithm of Section \ref{subsubs: VNS approach} is initialized with the greedy approach of Section \ref{subs: initialization procedure} to be then improved with the local search \emph{slack} from Section \ref{subsubs: slack local search}. Similar explanations have the rest of approaches, namely, \emph{VNS + greedy + sf} \emph{VNS + greedy + sr},  \emph{VNS + greedy + sc}, \emph{VNS + greedy + sfr} and \emph{VNS + greedy + sfrc}. The case with the genetic approach is also analogous. Note that, in order to get a better visualization of the results, we omit the boxplots for the methodologies \emph{naive}, \emph{VNS} and \emph{genetic}. Looking at Figures \ref{fig: boxplot 30c and 25p sensitivity local search VNS zoom} and \ref{fig: boxplot 30c and 25p sensitivity local search gen zoom}, it can be seen that in both cases, including the \emph{slack} local improvement increases the objective value compared with the corresponding approach without such improvement, i.e., compared with \emph{VNS + greedy} or \emph{gen + greedy}. Regarding the local searches that involve changes in the binary variables, \emph{fill}, \emph{reassignment} and \emph{conditional reassignment}, we observe that in Figure \ref{fig: boxplot 30c and 25p sensitivity local search VNS zoom} apply the \emph{fill} or the \emph{conditional reassignment} strategies after the \emph{slack} local searches have slightly better performance than applying the \emph{reassignment}. In other words, the \emph{VNS + greedy + sf} and \emph{VNS + greedy + sc} yield better objective values than \emph{VNS + greedy + sr}. In contrast, in the genetic case, the approach \emph{gen + greedy + sc} is better than \emph{gen + greedy + sf} and \emph{gen + greedy + sr}. Regarding the inclusion of the \emph{reassignment} local search to \emph{slack} and \emph{fill} in the VNS approach of Figure \ref{fig: boxplot 30c and 25p sensitivity local search VNS zoom}, i.e., analysing the strategy \emph{VNS + greedy + sfr}, we observe some improvement compared to the cases where only two local searches are used, namely the case \emph{VNS + greedy + sf}. In addition, slightly better results can be seen in the bottom whisker of the boxplot if we compare \emph{VNS + greedy + sfr} and \emph{VNS + greedy + sfrc}.  Similar conclusions are obtained when comparing \emph{gen + greedy + sfr} and \emph{gen + greedy + sfrc} in Figure \ref{fig: boxplot 30c and 25p sensitivity local search gen zoom}. More precisely, in this case, both approaches yield the same results.

Next, we analyze the evolution of different methodologies over time. In particular, we analyze the results of \emph{naive}, \emph{VNS + greedy + sfrc}, \emph{gen + greedy + sfrc}, and \emph{benchmark}. Moreover, we also include a variant of the \emph{exact} method, named \emph{exact\_h}, where \emph{h} comes from heuristics. The idea of this experiment is to compare our proposal, which is a heuristic approach, with the great amount of heuristic algorithms available in the literature and that are already included in Gurobi to help the resolution of MILPs. To this aim, in addition to tuning the parameter \texttt{Threads} to the value 1, we also modify the Gurobi parameter named \texttt{Heuristics} which controls the proportion of time spent in MIP heuristics. Such a parameter is set to its maximum value 1.

\begin{figure}[htb!]
\centering
\includegraphics[scale=0.7]{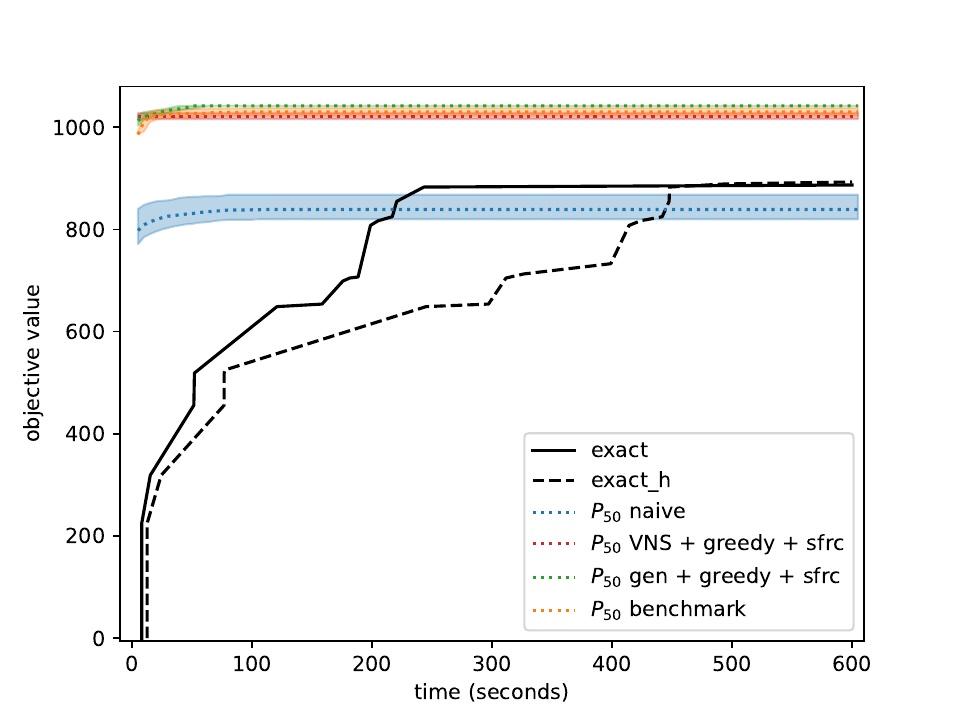}
\captionsetup{labelfont={color=black}}
\caption{Evolution in time of the results with the different approaches in the instance with 30 customers and 25 products.}
 \label{subfig: evolution max values 30c_25p}
    \end{figure}

\begin{figure}[htb!]
\centering
\includegraphics[scale=0.6]{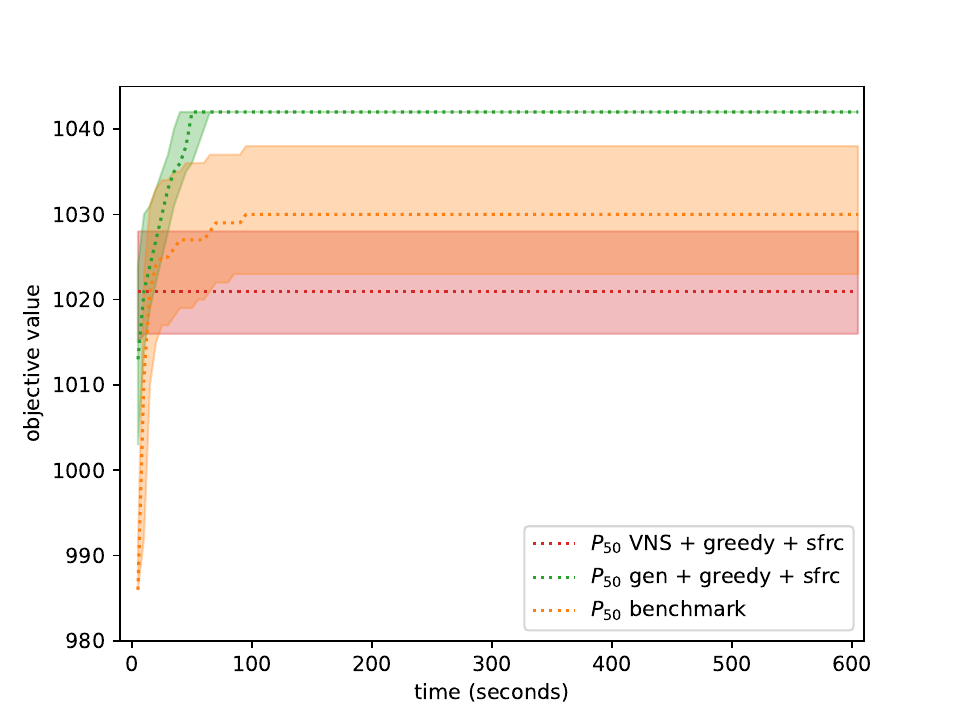}
\caption{Zoom of the evolution in time of the results with the different approaches in the instance with 30 customers and 25 products.}
 \label{fig: evolution max values 30c_25p zoom}
    \end{figure}

Figure \ref{subfig: evolution max values 30c_25p} depicts these results. The solid black line represents the evolution of the branch-and-bound algorithm carried out with the \emph{exact} approach, whereas the dashed black line indicates the behavior of the \emph{exact\_h} approach. In the first case, we can observe that the \emph{optimal} value of 887 is attained in approximately 200 seconds. However, the \emph{optimal} value of the exact\_h strategy coincides with 893 and is reached around 400 seconds. On the other hand, we represent the results of the\emph{naive}, \emph{VNS + greedy + sfrc}, \emph{gen + greedy + sfrc}, and \emph{benchmark} in blue, red, green and orange, respectively. In these experiments, we set the same stopping criterion as before, a maximum number of $24000$ points in the set $\mathcal{L}$. Note that the dotted lines, independent of which approach is used, indicate the median value (or the $50$-th percentile, $P_{50}$) obtained along the 1000 runs of the incumbent approach at each time instant. The shadowed areas represent the range where the $90\%$ of the results are located. That is to say, these areas are limited by the $5$-th percentile, $P_{5}$, and the $95$-th percentile, $P_{95}$.
We decided not to analyze the best and the worst cases yield along the time (which corresponds to $P_0$ and $P_{100}$, respectively) but the percentile $P_5$ and $P_{95}$. The reason for this choice is that we do not want to get wrong conclusions derived from the extreme values that may be the consequence of atypical results. This way, using $P_5$ and $P_{95}$ will provide more consistent conclusions. We can see that in the first seconds, all the approaches, \emph{naive}, \emph{VNS + greedy + sfrc}, \emph{gen + greedy + sfrc}, and \emph{benchmark} give better results than the \emph{exact} and \emph{exact\_h} strategies. Hence, if a feasible solution with a good objective value is desired, then the user can choose any of these methodologies. Besides, we observe that the shadowed areas associated with the \emph{VNS}-based approach (in red), the \emph{genetic}-based strategy (in green), and the \emph{benchmark} methodology (in orange) always take higher values than the ones of the \emph{naive} method (in blue). Indeed, the lower limit of the red, green or orange areas that corresponds to the $5$-th percentile of the \emph{VNS}-based, \emph{genetic}-based or \emph{benchmark} approaches, respectively, that is, to the (almost) worst case that we may obtain, has a value that, in the limit, are much larger to the one obtained by the \emph{exact} and \emph{exact\_h} methodologies. Moreover, note that, even though the three strategies, \emph{VNS + greedy + sfrc}, \emph{gen + greedy + sfrc} and \emph{benchmark} yield similar results in median, the variability is different among them. To better see this effect, we refer to Figure \ref{fig: evolution max values 30c_25p zoom}, where we observe that the \emph{gen + greedy + sfrc} approach is the one with lower variance. In particular, from the second 50 approximately, this methodology attains the value 1042 in all the simulations run.

Finally, we perform a sensitivity analysis of our proposal concerning the input parameter $L_0$, i.e., with respect to the initial number of vectors to be included in the set $\mathcal{L}$. See Algorithms \ref{alg: VNS} and \ref{alg: genetic} for more details. To this aim, we compare the approaches \emph{VNS + greedy + sfrc} and \emph{gen + greedy + sfrc} which uses $L_0= 1000$, with two new values, namely $L_0 = 100$ and $L_0 = 5000$. These strategies are also compared with the \emph{benchmark} method. The rest of the hyperparameters, different from $L_0$, remain unchanged. In other words, $Q = 100$ and $T = 500$ for the \emph{VNS + greedy + sfrc} algorithm, and $Q = 1000$ and $T=500$ for the \emph{gen + greedy + sfrc} approach, and for the \emph{benchmark} methodology, too. We remove the \emph{naive} strategy from the comparison since we observe that no matter which value is used to initialize $L_0$, the \emph{VNS}-based and the \emph{genetic}-based approaches got better results in terms of the objective values than those provided by \emph{naive}. Furthermore, the stopping criterion for all the approaches is again the cardinal of the set $\mathcal{L}$ obtained at the end of the algorithms, which is set to $24000$ points.

\begin{figure}[htb!]
\centering
\includegraphics[scale = 0.6]{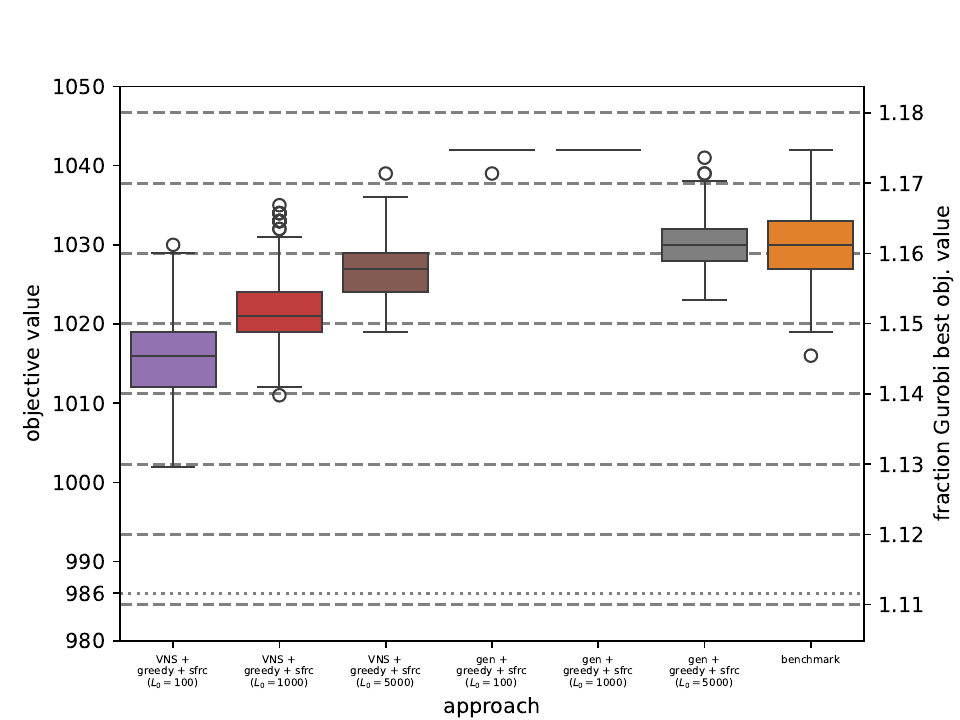}
\caption{Comparison of the approaches with different values of $L_0$ in the instance with 30 customers and 25 products.}
\label{fig: boxplot 30c and 25p comparison zoom}
\end{figure}

Figure \ref{fig: boxplot 30c and 25p comparison zoom} shows the boxplots of objective values $\phi^{max}$ obtained in the 1000 runs of the different approaches. The objective value provided by Gurobi using the \emph{exact} algorithm is 887. The objective value provided by the first point of the greedy initialization is 986, which is marked with a horizontal dotted line. The interpretation of the two vertical axes is exactly the same as the one in Figure \ref{subfig: boxplot 30c and 25p}. Firstly, we observe that no matter which value is used to initialize $L_0$, the \emph{VNS}-based and \emph{genetic}-based approaches got better results in terms of the objective values than those provided by Gurobi. The same behavior can be seen with the \emph{benchmark} strategy. Besides, we can see in the \emph{VNS}-based case that the larger the initialization value, the higher the central box of the boxplots is located. In other words, higher values of $L_0$ have bigger objective values in the 50\% central of the observations. This phenomenon may be due to the fact that, when increasing the value $L_0$, we sample more vectors of prices whose associated upper-level objective functions may have larger values. Therefore, the initial set $\mathcal{Q}$ may contain vectors with higher objective values, implying that the VNS approach is initialized with better local optima and, hence, gets better results. 

A different situation occurs in the \emph{genetic}-based simulations. It can be observed that the larger $L_0$, the smaller the objective values of the upper-level function. Particularly, this difference is more significant if we compare the results of $L_0 = 100$ and $L_0 = 1000$ with those of $L_0 = 5000$. The reason for this phenomenon may be due to the evolutionary process that the genetic algorithms follow. Sampling more initial points ($L_0 = 5000$) implies to have more vectors to choose from in the crossover and mutation operators, and therefore, one may have more difficulties in finding good objective values. In addition, we see that changing from $L_0 = 100$ to $L_0 = 1000$ no have strong consequences in terms of the objective values. Only in few instances of the case $L_0 = 100$, the algorithm is not able to attain the true optima 1042.

\subsection{Case study with 60 customers and 50 products}

This instance of the Rank Pricing Problem is formed by $K = 60$ customers and $I = 50$ products. We compare our approaches with the \emph{benchmark} strategy as well as with the \emph{exact} method. We first set as stopping criteria the cardinal of the set $\mathcal{L}$ which is fixed to 24000. The rest of the input parameters of our proposal are set to their default values, that is to say, $L_0 = 1000$, $Q = 100$, and $T = 500$ in the \emph{VNS}-based strategies and $L_0 = 1000$, $Q = 1000$, and $T = 500$ for the \emph{genetic}-based approaches. The \emph{optimal} objective value obtained by the \emph{exact}  approach is equal to 1151.

\begin{figure}[htb!]
\centering
\includegraphics[scale = 0.6]{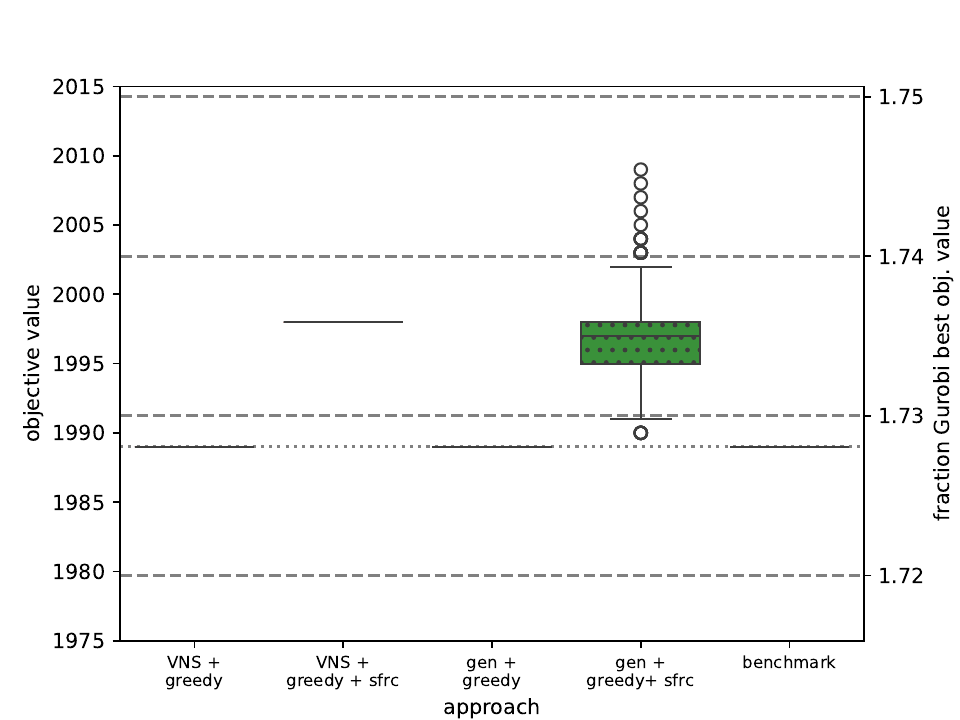}
\caption{Comparison of the different approaches for the instance with 60 customers and 50 products using as stopping criteria that the cardinal of the set $\mathcal{L}$ reaches 24000 points.}
\label{fig: boxplot 60c and 50p 24000p zoom}
\end{figure}

\begin{figure}[htb!]
\centering
\includegraphics[scale = 0.7]{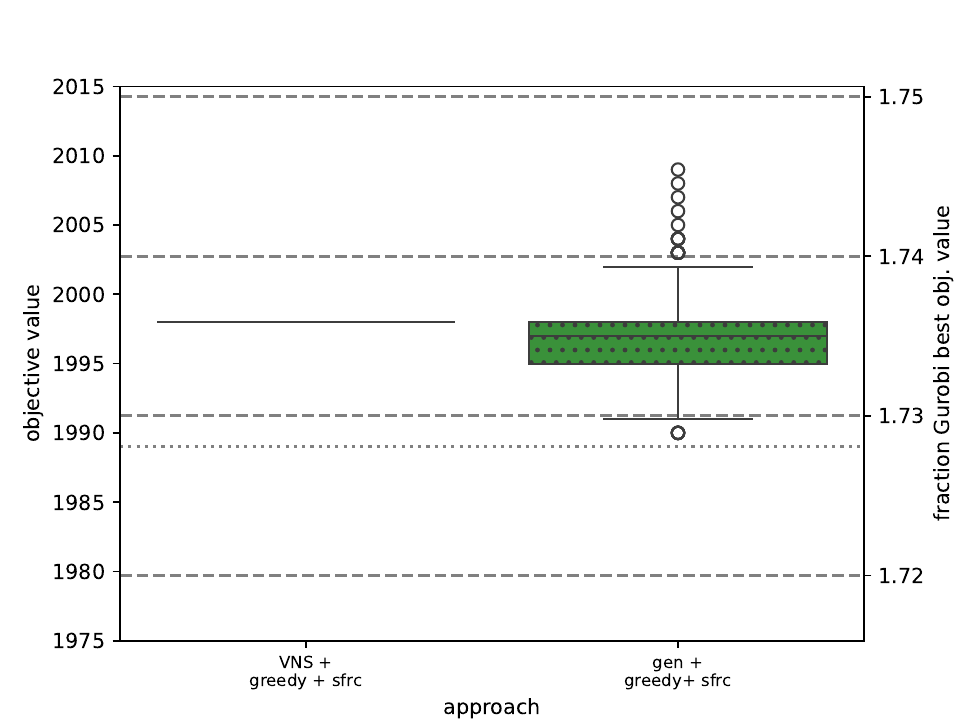}
\captionsetup{labelfont={color=black}}
\caption{Zoom of the comparison of the different approaches for the instance with 60 customers and 50 products using as stopping criteria that the cardinal of the set $\mathcal{L}$ reaches 24000 points.}
\label{subfig: boxplot 60c and 50p 24000p zoom}
    \end{figure}

Figure \ref{fig: boxplot 60c and 50p 24000p zoom} shows the results of the \emph{naive}, \emph{VNS + greedy}, \emph{VNS + greedy + sfrc}, \emph{gen + greedy}, \emph{gen + greedy + sfrc} and \emph{benchmark} approaches. For the sake of visualization, we remove the results from \emph{naive}, \emph{VNS} and \emph{genetic} approaches. We observe that all the strategies get larger objective values than the one provided by Gurobi using the single-level reformulation \eqref{eq: single-level reformulation}. Regarding the effect of the greedy initialization in the \emph{VNS}-based or \emph{genetic}-based approaches, we observe that the use of the greedy approach has indeed a positive effect in terms of the objective values. In particular, both methodologies always attain the value 1989 (plotted with a dotted line) which is the objective value obtained with the first point of the \emph{greedy} initialization. This means that neither applying the \emph{VNS} or the \emph{genetic} strategies will imply improvements in the values. A similar situation occurs in the \emph{benchmark} approach, where apart from the \emph{genetic} with the \emph{greedy} initialization, the \emph{optimization-based} local search of Section \ref{subsubs: optimization-based local search} is run. It means that such local search provides no improvement regarding the objective value of the initialization. Let us now focus on the effect of the proposed local searches on the base heuristics. In both cases, \emph{VNS + greedy + sfrc} and \emph{gen + greedy + sfrc} have better values than their corresponding approaches without the local searches, then demonstrating the effectiveness of the local improvements proposed in this paper. Let us focus with more details on the differences between these two approaches. To this aim, let us focus on Figure \ref{subfig: boxplot 60c and 50p 24000p zoom}, where a zoom of these values is made. It can be seen that, in this case, the \emph{VNS}-based setting yields results with a smaller variance than those of the \emph{genetic}-based methodology. In fact, the 1000 simulations of the \emph{VNS} approach are equal to the value 1998, so the variance in this methodology is zero. In contrast, some runs of the \emph{genetic} algorithm yields larger values than those provided by the VNS, i.e., larger than 1998, even reaching objective values close to 2010. In addition, approximately the $75\%$ of the objective values obtained by the \emph{genetic}-based methodology are smaller than the value obtained by the VNS counterpart (1998). Whether to apply the \emph{VNS} or the \emph{genetic} strategies combined with the \emph{greedy} initialization and the local searches \emph{sfrc} is an user decision. If one wants to have some chance of getting a value larger than 1998, then the \emph{genetic} algorithm should be used. In contrast, if the user wants to have small variance in the results, the \emph{VNS} setting is to be applied. 

In what follows, we will compare the results obtained with two different stopping criteria, namely the cardinal of 24000 points for the set $\mathcal{L}$, and the time limit of $600$ seconds. For the sake of comparison, some remarks should be done before continuing. Firstly, we want to highlight that, the computational time of all algorithms when the stopping criterion of 24000 points is used is around  100 seconds. Secondly, note that we have chosen to use also the time limit of $600$ seconds because this is the computational time allowed in Gurobi for solving the \emph{exact}approach. Finally, note that in the \emph{VNS}-based, \emph{genetic}-based and \emph{benchmark} methodologies, we measure the computational time just taking into account the time provided by Gurobi when solving the lower-level problems \eqref{eq: lower level}. We are aware that there are different parts of the algorithm that require some non-significant time, such as, for instance, the local searches. However, we think that the measure of such a computational time depends on the way the approach is coded. In contrast, the time provided by Gurobi only depends on the solver itself. Consequently, fair conclusions are derived.

\begin{figure}[htb!]
\centering
\includegraphics[scale = 0.6]{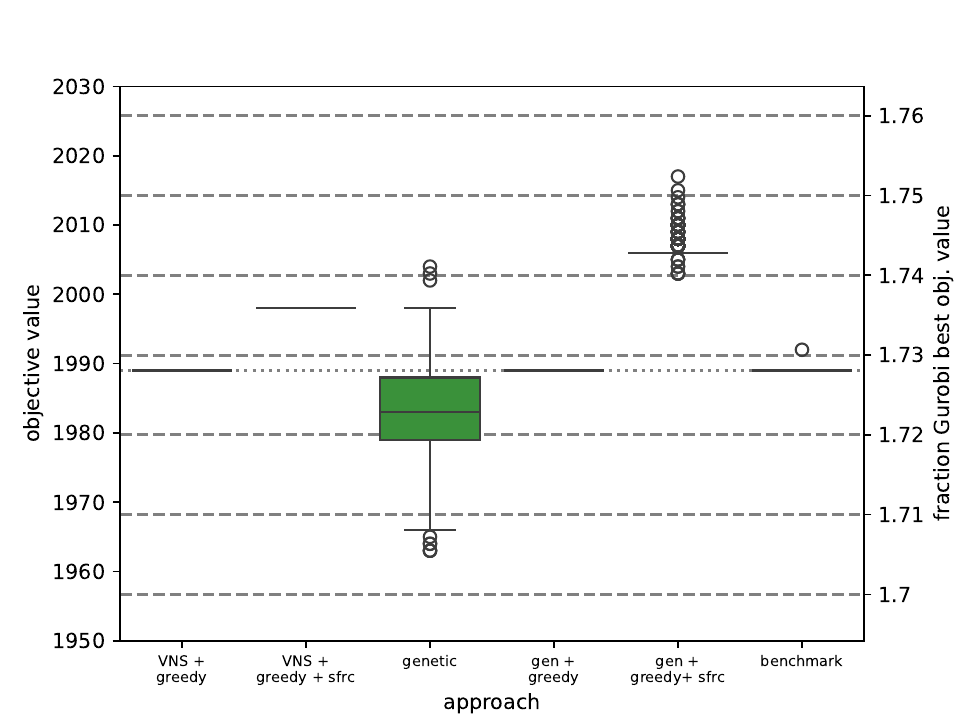}
\caption{Comparison of the different approaches for the instance with 60 customers and 50 products using as stopping criteria a time limit of 600 seconds.}
\label{fig: boxplot 60c and 50p 600s zoom}
\end{figure}

Figure \ref{fig: boxplot 60c and 50p 600s zoom} shows the results obtained for the different approaches (except \emph{naive} and \emph{VNS}) using the computational time limit of 600 seconds as stopping criterion. As expected, the objective values of Figure \ref{fig: boxplot 60c and 50p 600s zoom} are larger (or in some cases equal) than those obtained in Figure \ref{fig: boxplot 60c and 50p 24000p zoom}, where we stopped the algorithm when the set $\mathcal{L}$ has 24000 points. For instance, we observe in Figure \ref{fig: boxplot 60c and 50p 600s zoom} that all the solutions of the methodology \emph{gen + greedy + sfrc} has an objective value greater than 2000, whereas in the case of Figure \ref{fig: boxplot 60c and 50p 24000p zoom} more than $75\%$ of the results obtained are below the value 2000. Moreover, when the \emph{VNS}-based strategies are used, it can be seen that there are no significant modifications when changing the stopping criterion. In particular, the results of \emph{VNS + greedy} and \emph{VNS + greedy + sfrc} are exactly the same in both figures. In contrast, some differences can be observed in the case of the \emph{genetic}-based and \emph{benchmark} strategies. Finally, the \emph{gen + greedy} method do not change their value, 1989, when the stopping criteria is modified.

\section{Conclusions and Further Research} \label{sec: conclusions}

This paper addresses the Rank Pricing Problem (RPP), a challenging NP-hard combinatorial optimization problem formulated as a bilevel model with lower-level binary variables. To tackle the RPP, we propose a heuristic approach that effectively leverages the problem's structure. Our method combines a well-established heuristic algorithm—specifically, we test standard implementations of a genetic algorithm and a Variable Neighborhood Search routine—with four novel local search subroutines designed to enhance the solutions provided by the base heuristic. Unlike previous approaches, our method avoids solving additional optimization problems, significantly reducing computational complexity.  

We evaluate our approach on various RPP instances with different numbers of products and clients. Results demonstrate that our method consistently outperforms MIP solvers within a short computational time. For example, in large instances, the most competitive variants of our strategy achieve objective values that are systematically more than 1.7 times higher than those obtained by Gurobi in under 600 seconds.  

Future research could explore further enhancements, such as integrating our heuristic approach with other optimization techniques to tackle even more complex bilevel problems or incorporating machine learning to refine decision-making in combinatorial models. Additionally, our method is well-suited for parallel computing. Another promising avenue for future work is adapting our approach to variants of the RPP where clients' budgets and/or preferences are unknown.

\section*{Acknowledgements}

This work was supported in part by the Spanish Ministry of Science and Innovation (AEI/10.13039/501100011033) through project PID2023-148291NB-I00. The authors thankfully acknowledge the computer resources (Picasso Supercomputer), technical expertise, and assistance provided by the SCBI (Supercomputing and Bioinformatics) center of the University of Málaga.

\bibliographystyle{apalike}
\bibliography{Bibliography_all}

\end{document}